%% file: main.tex
\documentclass[11pt]{amsart} 
\usepackage[utf8]{inputenc} 
\usepackage{amssymb, amsmath}
\usepackage{bbm}
\usepackage{dsfont}
\usepackage{amsthm}
\usepackage{amsrefs} 
\usepackage{graphicx} 
\usepackage{enumitem} 
\usepackage{subfigure}
\usepackage{fancyvrb, verbatim} 
\usepackage{graphicx}
\usepackage{subcaption}
\usepackage{hyperref} 
\usepackage{cleveref,cite,url}

\hypersetup{
    colorlinks=true,
    linkcolor=cyan,
    urlcolor=cyan,
    citecolor=cyan,
    }
    
\usepackage{cleveref,cite,url}
\usepackage{hologo} 
\usepackage{comment} 
\usepackage{hologo}
\usepackage{comment} 
\usepackage{xcolor}

\newcommand{\CC}{\mathbb{C}} 
 
\DeclareMathOperator{\rank}{rank}

\theoremstyle{plain} 
\newtheorem{thm}{Theorem}[section] 
\newtheorem*{thm:infinitelymanyprimes}{Theorem \ref{thm:infinitelymanyprimes}} 
\newtheorem{lem}{Lemma}[section]

\newtheorem{cor}{Corollary}[section]

\theoremstyle{definition}
 
\newtheorem{ex}[thm]
{Example}
 
\theoremstyle{remark}
\newtheorem{rmk}{Remark}[section]

\input{notation}

\begin{document}

\title{On the exactness of Non-unique game SDP over finite groups}

\author{Diego Cifuentes}
\address{~}
\email{diegcif@gmail.com}

\author{Zhuorui Li}
\address{Georgia Institute of Technology \ Atlanta, GA, USA}
\email{zli966@gatech.edu}

\date{}

\begin{abstract}
We study the exactness of semidefinite programming (SDP) relaxations for non-unique games (NUGs) over finite groups. We first provide a characterization of the exact SDP solutions and derive a deterministic sufficient condition for  exactness. We then apply this condition to multireference alignment and joint assignment, showing that the SDP achieves minimax-optimal thresholds for exact recovery under the respective noise models. We also investigate an iterative rounding procedure that promotes SDP solutions of minimal rank. Numerical experiments on Max-2-Lin, joint assignment and multireference alignment, demonstrate the effectiveness of our rounding approach.
\end{abstract}

\maketitle

\section{introduction}

The \emph{non-unique games} (NUG) framework~\cite{bandeira2020non} provides a unified formulation for a broad class of optimization problems arising in combinatorial optimization, machine learning, statistical physics, signal processing, and computer vision, which include \texttt{Max-2-Lin}~\cite{andersson2001new,khot2007optimal,feige2004systems}, the Ising and Potts models~\cite{barahona1982computational,wu1982potts}, and certain metric-labeling and maximum a posteriori inference problems in graphical models~\cite{kleinberg2002approximation,boykov2001fast,kolmogorov2004energy,wainwright2005map}. Other important applications include multireference alignment~\cite{bandeira2014multireference,bendory2017bispectrum,perry2019sample}, registration of point clouds~\cite{chaudhury2015global}, joint assignment from pairwise differences~\cite{chen2018projected}, angular synchronization~\cite{singer2011angular}, and orientation estimation in cryo-electron microscopy (cryo-EM)~\cite{singer2011three,bandeira2020non,shkolnisky2012viewing,wang2013exact,wang2013orientation}.

In this paper, we study the non-unique games (NUG)
problem over a finite group $G$.
Let $e$ denote the identity element of $G$,
and let $g_i\in G$, $i=1,\ldots,n$, be unknown
group elements.
For each ordered pair $(i,j)$, suppose we are
given a pairwise potential
\[
f_{ij}: G\to\mathbb R.
\]
The general problem is
\begin{equation}
\tag{$\mathrm{P}_{\RR}$}
\min_{ g\in \cG^n}
\sum_{i,j=1}^n f_{ij}(g_i g_j^{-1}), 
\label{eq:cyclic-difference}
\end{equation}
where $ g=(g_1,\ldots,g_n)$ and
$g_i g_j^{-1}$ denotes the relative group element
between $g_i$ and $g_j$.

This formulation captures models in which the interaction between two variables depends only on their relative group element. Since the objective depends only on $g_i g_j^{-1}$, it is invariant under a common right multiplication of all labels. 
To understand \eqref{eq:cyclic-difference} better, we present a concrete instance of \eqref{eq:cyclic-difference}. 

One important special case is \texttt{Max-2-Lin} over the cyclic group
\[
\mathbb Z_L:=\{0,1,\dots,L-1\}.
\]
In this problem, we are given weighted equations
\[
s_i-s_j\equiv c_{ij}\bmod L,
\qquad (i,j)\in E,
\]
with weights \(w_{ij}\ge 0\). The goal is to maximize the total weight of satisfied equations:
\begin{equation}
\max_{s\in\mathbb Z_L^n}
\sum_{(i,j)\in E}
w_{ij}\mathds{1}\{s_i-s_j\equiv c_{ij}\bmod L\}.
\label{eq:max2lin}
\end{equation}
This fits into \eqref{eq:cyclic-difference} by defining 
\[
f_{ij}(t)=-w_{ij}\mathds{1}\{t=c_{ij}\}.
\]
Then maximizing \eqref{eq:max2lin} is equivalent to minimizing
\[
\sum_{(i,j)\in E} f_{ij}(s_i-s_j).
\]
Thus, \texttt{Max-2-Lin} is a special case of \eqref{eq:cyclic-difference}.

However, solving \eqref{eq:cyclic-difference} is computationally intractable in general. To address this challenge, Bandeira et al.~\cite{bandeira2020non} developed a semidefinite programming relaxation framework (NUG-SDP) for non-unique games. Despite its broad applicability, conditions under which the NUG-SDP exactly solves \eqref{eq:cyclic-difference} remain largely unexplored. In this paper, we provide a deterministic sufficient condition that guarantees the exactness of the NUG-SDP over finite groups. 
When $G$ is an infinite compact group, such as $SO(2)$ and $SO(3)$, NUG-SDP can still be applied by truncating the Fourier expansion to finitely many irreducible representations. However, although the truncated SDP has been demonstrated to be effective in practice \cite{bandeira2020non}, its optimal solution is not necessarily an optimal solution of the original problem \eqref{eq:cyclic-difference}. Building on our exactness condition, we further show that the NUG-SDP achieves minimax-optimal exact recovery for two important applications: multireference alignment and joint assignment. Finally, when the SDP solution does not directly yield a group assignment, we propose a rounding method and demonstrate its
effectiveness through numerical experiments.

\subsection*{Related works}

Closely related works arise in angular and phase synchronization, where unknown phases are estimated from noisy measurements of their pairwise ratios. Singer \cite{singer2011angular} proposed spectral and semidefinite methods for angular synchronization, and subsequent work established conditions under which the maximum-likelihood SDP relaxation is tight \cite{bandeira2017estimation}. Gao and Zhang \cite{gao2022sdp} proved the minimax optimality of an SDP relaxation for phase synchronization. Beyond single complex quadratic programs, Gao and Zhao \cite{gao2019multifrequency} introduced a multi-frequency phase synchronization framework that combines measurements from several frequency channels. However, the above works mainly focus on optimization over
continuous phase variables, whereas we study non-unique games
over finite groups.

Motivated by problems in signal processing, physics, and theoretical computer science, many works study single frequency complex quadratic programs over roots of unity \cite{ma2004semidefinite,de2009design,arora2021efficient,mobasher2007near,soltanalian2014designing,sinjorgo2024cuts}. Semidefinite relaxations and randomized rounding schemes for this setting were developed in \cite{zhang2006complex,so2007approximating}, with particular emphasis on approximation guarantees for a single complex quadratic form. Lu et al. \cite{lu2019tightness} proposed an SDP to solve an inhomogeneous complex quadratic program over roots-of-unity, and explored its tightness for MIMO model. Al-Sulami et al. \cite{al2026convergence} studied sum-of-square hierarchy to solve such a problem and provided convergence results for the hierarchy. Later, Al-Sulami et al. \cite{al2026zonotope} also studied a binary integer programming reformulation to find the exact solution of this problem. In a related direction, Fawzi et al. \cite{fawzi2016sparse} studied sparse sum-of-squares certificates for nonnegative functions over finite abelian groups.

Computationally efficient nonconvex methods have also been studied as alternatives to SDP relaxations~\cite{boumal2016nonconvex,liu2017estimation,bendory2017bispectrum,lerman2022robust}. Some of these methods have strong theoretical guarantees and achieve minimax optimality in specific applications. Chen and Cand\`es~\cite{chen2018projected} proposed a projected power method for joint alignment from pairwise differences, a special case of \eqref{eq:cyclic-difference}, and established information-theoretically optimal exact-recovery guarantees. More recently, Gao and Zhang~\cite{gao2022iterative} developed a general iterative framework for discrete structure recovery, encompassing multireference alignment and group synchronization, and proved that their algorithms achieve minimax-optimal error rates under suitable conditions.

The work most closely related to ours is that of Bandeira et al.~\cite{bandeira2020non}, who generalized their earlier SDP relaxation for multireference alignment~\cite{bandeira2014multireference} to the non-unique games framework over compact groups. Motivated by applications involving special orthogonal groups \cite{shkolnisky2012viewing,singer2011three,wang2013exact,saunderson2015semidefinite}, they developed an ADMM-based solver and demonstrated its robustness in low SNR problems over $SO(2)$ and $SO(3)$ in comparison with classic synchronization approaches.  More recently, Shi et al.~\cite{shi2026orientation} extended the NUG framework to accommodate symmetries in cryo-EM orientation estimation and proposed a numerical refinement procedure. In contrast, our work focuses on establishing a deterministic sufficient condition for the exactness of the NUG-SDP when $G$ is finite and, building on this condition, deriving minimax-optimal exact-recovery guarantees for multireference alignment and joint assignment. When the SDP relaxation is not tight, we also propose a rounding method  and demonstrate its effectiveness in numerical experiments.

\subsection*{Contribution and organization}

In this paper, we study the NUG-SDP relaxation for \eqref{eq:cyclic-difference}. Our main contributions are summarized as follows:

\begin{itemize}
\item \textit{A deterministic exactness condition}. In \Cref{sec:exactness}, we analyze the NUG-SDP relaxation and provide a deterministic condition which guarantees the exactness of this SDP.

\item \textit{Multireference alignment}. In \Cref{sec:mra}, we apply our condition to multireference alignment and show that minimax optimality is achieved by NUG-SDP. 

\item \textit{Joint assignment}. In \Cref{sec:joint_assignment}, we specialize our exactness condition to the joint assignment problem and show that minimax optimality is achieved by NUG-SDP. 

\item \textit{Iterative rounding}. In \Cref{sec:rounding} we propose an iterative rounding method for NUG-SDP when the relaxation is not exact, and show that every minimal rank feasible solution is a fixed point of the iteration.

\item \textit{Numerical experiments}. In \Cref{sec:experiments}, we conduct numerical experiments on \texttt{Max-2-Lin}, joint assignment, and multireference alignment to evaluate the performance of NUG-SDP and our rounding method.
\end{itemize}

\subsection*{Notation}

Let $\CC^n$ be the space of complex vectors of length~$n$, $\mathbb{H}^n$ the space of $n \times n$ Hermitian matrices, and $\mathbb{H}^n_+$ the set of $n \times n$ PSD matrices. Denote $\omega = e^{\frac{2\pi \boldsymbol{\mathrm{i}}}{L}}$ as the $L$-th root of unity. Denote $[L] = \{1,2,\dots,L\}$. Let \(\vone_n\) be an all-ones vector of length~$n$, and \(E_n = \vone_n\vone_n^T\). Given \(z\in\CC^n\), let $\overline{z}$ be its complex conjugate. Denote by $A\otimes B$ as the Kronecker product of $A$ and $B$. Given a matrix $Z$, $\diag(Z)$ extracts the diagonal of $Z$ into a vector, and $\offdiag(Z)$ keeps the off-diagonal entries of $Z$ and sets the diagonal entries to $0$. Given a matrix $Z \in \CC^{m\times n}$, let $\|Z\|$ denote its spectral norm;  let $Z^T$ and $Z^H$ denote its transpose and conjugate transpose, respectively.

\section{SDP relaxation}\label{sec:SDP}

In this section, we review the NUG-SDP relaxation as proposed in \cite{bandeira2020non} over a
general finite group $G$. We first reformulate
\eqref{eq:cyclic-difference} as an equivalent quadratic
program using the irreducible unitary representations of
$G$. We then introduce the corresponding NUG-SDP
relaxation.

\subsection{An equivalent reformulation via group representations}

Given potential functions $f_{ij}(\cdot)$, for $i,j\in[n]$,
define matrices $F_g\in\RR^{n\times n}$, $g\in G$, by
$(F_g)_{ij}=f_{ij}(g)$. Without loss of generality, we assume
$f_{ij}(g)=f_{ji}(g^{-1})$ for all $i,j\in[n]$ and $g\in G$.
Then \eqref{eq:cyclic-difference} can be written as
\begin{equation}
\min_{g\in G^n}
\sum_{i,j=1}^n (F_{g_i g_j^{-1}})_{ij}.
\label{eq:F-form}
\end{equation}

Let $\widehat G$ denote a complete set of inequivalent
irreducible unitary representations of $G$. For each $\pi\in\widehat G$, let
$\pi:G\to U(d_\pi)$ be a representation of dimension $d_\pi$,
where $U(d_\pi)$ denotes the group of $d_\pi\times d_\pi$
unitary matrices. Let $\pi_0$ denote the trivial representation. See \cite{steinberg2012representation,serre1977linear} for more details.

Applying the Fourier transform $\cF_G$ over $G$ to the pairwise
potentials, define matrices
\begin{equation}\label{eq:group_fourier}
\C:=(C_\pi)_{\pi\in\widehat G}
=\mathcal F_G((F_g)_{g\in G}),
\quad \text{where}~
C_\pi:=\sum_{g\in G}F_g\otimes\pi(g).
\end{equation}
Then, by Fourier inversion $\cF_G^{-1}$, we have
\[
(F_g)_{g\in G} = \cF_G^{-1}\big((C_\pi)_{\pi\in\widehat G}\big)
\]
where
\begin{equation}
(F_g)_{ij}
=
\frac{1}{|G|}
\sum_{\pi\in\widehat G}
d_\pi\operatorname{Tr}
\left(\pi(g)^H(C_\pi)_{ij}\right),
\qquad g\in G.
\label{eq:F-inverse}
\end{equation}

When $G = \ZZ_L$, we have $\pi_\ell(k) = \omega^{-k\ell}$, and the Fourier transform $\cF_{\ZZ_L}$ is the standard discrete Fourier transform, that is, 
\begin{equation*}
\C := (C_0,\dots,C_{L-1}) = \cF_{\ZZ_L}(F_0,\dots,F_{L-1}) 
\end{equation*}
where \(C_\ell = \sum_{k=0}^{L-1}\omega^{-k\ell}F_k\). And the inverse Fourier transform $\cF_{\ZZ_L}^{-1}$ is
\begin{equation*}
\boldsymbol{F} := (F_0,\dots,F_{L-1}) = \cF^{-1}_{\ZZ_L}(C_0,\dots,C_{L-1}) 
\end{equation*}
where \(F_\ell = \frac{1}{L}\sum_{k=0}^{L-1}\omega^{k\ell}C_k\).
% We discuss applications over the cyclic group $\ZZ_L$ in \Cref{sec:mra} and \Cref{sec:joint_assignment}.
For each $\pi\in\widehat G$, let
\[
V_\pi(g)
=
\begin{pmatrix}
\pi(g_1)\\
\vdots\\
\pi(g_n)
\end{pmatrix}
\in\mathbb C^{nd_\pi\times d_\pi}.
\]
Substituting $g=g_i g_j^{-1}$ into \eqref{eq:F-inverse},
we obtain
\[
\begin{aligned}
(F_{g_i g_j^{-1}})_{ij}
&=
\frac{1}{|G|}
\sum_{\pi\in\widehat G}
d_\pi\operatorname{Tr}
\left(
\pi(g_i)^H(C_\pi)_{ij}\pi(g_j)
\right).
\end{aligned}
\]
Therefore, \eqref{eq:F-form} is equivalent to
\begin{equation}\tag{$\mathrm{P}_{\CC}$}
\min_{g\in G^n}
\sum_{\pi\in\widehat G\setminus\{\pi_0\}}
d_\pi\operatorname{Tr}
\left(V_\pi(g)^H C_\pi V_\pi(g)\right),
\label{eq:multi-frequency-nonconstant}
\end{equation}
where the trivial representation $\pi_0$ is excluded,
since its contribution to the objective is constant,
and the positive scaling factor $1/|G|$ is omitted.

\subsection{Convex relaxation}

In this section, we present the non-unique games (NUG) SDP
as proposed in \cite{bandeira2020non}. The main obstacle
to solving \eqref{eq:multi-frequency-nonconstant} is the
following nonconvex set
\[
\cQ_G^n :=
\left\{
\Z=(Z_\pi)_{\pi\ne\pi_0}:
Z_\pi=V_\pi(g)V_\pi(g)^H,\ 
\forall \pi\ne\pi_0,\ g\in G^n
\right\}.
\]

The key idea is to derive a tractable convex relaxation
of $\cQ_G^n$. To this end, we define the inverse Fourier
transform of $Z$ by
\[
(\hat Z_g)_{ij}
=
\frac{1}{|G|}
\sum_{\pi\in\widehat G}
d_\pi\operatorname{Tr}
\left(\pi(g)^H(Z_\pi)_{ij}\right),
\qquad g\in G,
\]
where $Z_{\pi_0}=E_n:=\vone_n\vone_n^\top$.
The NUG-SDP over $G$ relaxes $\cQ_G^n$ using the
convex set $\cZ_G^n$ which is defined as follows:
\[
\begin{aligned}
\cZ_G^n:=\Big\{
\Z=(Z_\pi)_{\pi\ne\pi_0}:\quad
&(\hat Z_g)_{g\in G}
=\mathcal F_G^{-1}
(E_n,(Z_\pi)_{\pi\ne\pi_0}),\\
&\diag(\hat Z_e)=\vone_n,\quad
\diag(\hat Z_g)=\vzero_n,\quad g\ne e,\\
&\hat Z_g=\hat Z_{g^{-1}}^T,\quad g\in G,\\
&\hat Z_g\in\RR_+^{n\times n},\quad g\in G
\Big\}.
\end{aligned}
\]
Note that $\hat Z_g=\hat Z_{g^{-1}}^T$ implies that
$Z_\pi$ is Hermitian for every $\pi\in\widehat G$.

Let $\C\bullet \Z:=
\sum_{\pi\in\widehat G\setminus\{\pi_0\}}
d_\pi\langle C_\pi,Z_\pi\rangle$.
Based on the above convex set $\cZ_G^n$, it is natural to consider the following SDP relaxation:
\begin{equation}\label{eq:sdp}\tag{Q}
\begin{aligned}
\min_{Z}\quad&
\C\bullet \Z\\
\text{s.t.}\quad&
Z\in\cZ_G^n,\\
&Z_\pi\succeq0,\quad
\pi\in\widehat G\setminus\{\pi_0\},
\end{aligned}
\end{equation}

The above SDP involves one PSD constraint for each
nontrivial irreducible representation $\pi\in\widehat G$,
with matrix dimension $nd_\pi$.
When $G$ is finite, only \((|\widehat G|-1+r)/2\) constraints are needed in implementation, where $r$ is the number of nontrivial self-conjugate irreducible representations of $G$; see \Cref{append:less_PSD}.

The relaxation is \textit{exact} or \textit{tight}
if an optimal solution of \eqref{eq:sdp}
belongs to the original feasible set $\cQ_G^n$. In the next section, we provide a deterministic
sufficient condition under which the SDP relaxation
\eqref{eq:sdp} exactly solves the original problem
\eqref{eq:multi-frequency-nonconstant}.

\section{deterministic exactness condition}\label{sec:exactness}

In this section, we investigate the condition under which the SDP relaxation \eqref{eq:sdp} is tight. First, we show that the relaxation is tight if and only if the solution has minimal rank. Then we further give a deterministic condition on the exactness of the SDP. 

\subsection{Exactness condition}

Given $\Z=(Z_\pi)_{\pi\ne\pi_0}\in\cZ_G^n$,
we say that $\Z$ is of \emph{minimal rank} if
$\rank(Z_\pi)=d_\pi$ for every
$\pi\in\widehat G\setminus\{\pi_0\}$. For $Z_\pi\in\mathbb H^{nd_\pi}$, let $(Z_\pi)_{ij}$ denote
its $(i,j)$-th $d_\pi\times d_\pi$ block.
We show that this condition is equivalent to
$Z$ belonging to the original feasible set $\cQ_G^n$.

\begin{thm}\label{thm:rank-one}
For $\Z=(Z_\pi)_{\pi\ne\pi_0}\in\cZ_G^n$, we have
\[
\Z\in\cQ_G^n\quad\Longleftrightarrow\quad
\rank(Z_\pi)=d_\pi,\quad
\forall\pi\in\widehat G\setminus\{\pi_0\}.
\]
\end{thm}

\begin{proof}
($\Leftarrow$) Suppose that $\rank(Z_\pi)=d_\pi$ for every
$\pi\in\widehat G\setminus\{\pi_0\}$. Let 
\[(\hat Z_g)_{g\in G} =
\mathcal F_G^{-1}(E_n,(Z_\pi)_{\pi\ne\pi_0}).
\]

We first show that $\hat Z_g\in\{0,1\}^{n\times n}$ for every $g\in G$. From the definition of $\cZ_G^n$, it follows that  for all $i\in[n]$, 
\[
(Z_\pi)_{ii}= \sum_{g\in G} (\hat Z_g)_{ii} \pi(g)= (\hat Z_e)_{ii} \pi(e) = I_{d_\pi}
\]
Suppose that
\[
Z_\pi = \begin{pmatrix}
    I_{d_\pi} & B^H\\
    B & D
\end{pmatrix}
\]
By the Schur complement rank identity, we have 
\[
\rank(Z_\pi) = d_\pi + \rank(D - BB^H)
\]
Since $Z_\pi$ has rank $d_\pi$, it admits a factorization
\[
Z_\pi=V_\pi V_\pi^H,
\qquad
V_\pi=
\begin{pmatrix}
V_{\pi,1}\\
\vdots\\
V_{\pi,n}
\end{pmatrix},
\]
where $V_{\pi,i}\in U(d_\pi)$ for every $i\in[n]$.
Consequently, we have
\[
(Z_\pi)_{pq}=V_{\pi,p}V_{\pi,q}^H,
\]
which implies
\begin{equation*}\label{eq:rank_Z_i}
\begin{aligned}
\|(Z_\pi)_{pq}\|_F^2 = \|V_{\pi,p} V_{\pi,q}^H\|_F^2 & =  \tr(V_{\pi,p} V_{\pi,q}^HV_{\pi,q} V_{\pi,p}^H) \\
& = \tr(V_{\pi,p}V_{\pi,p}^H)\\
& = \tr(I_{d_\pi}) = d_\pi,
\qquad \forall p,q\in[n].
\end{aligned}
\end{equation*}

For a fixed pair $(p,q)$, let
$\lambda_{pq}(g):=(\hat Z_g)_{pq}\geq0$.
By the definition of the Fourier transform,
\[
(Z_\pi)_{pq}
=
\sum_{g\in G}\lambda_{pq}(g)\pi(g).
\]
For the trivial representation $\pi_0$, we have
\[
\sum_{g\in G}\lambda_{pq}(g)
=(Z_{\pi_0})_{pq}=1.
\]
Applying Plancherel's formula over finite groups\footnote{More details are introduced in \Cref{append:Plancherel}},
we obtain
\[
\begin{aligned}
\sum_{g\in G} (\hat Z_g)^2_{pq}=\sum_{g\in G}\lambda^2_{pq}(g)
&=
\frac{1}{|G|}
\sum_{\pi\in\widehat G}
d_\pi\|(Z_\pi)_{pq}\|_F^2\\
&=
\frac{1}{|G|}
\sum_{\pi\in\widehat G}d_\pi^2
=1.
\end{aligned}
\]
Since we have
$\sum_{g\in G}\lambda_{pq}(g)=1$, and $\sum_{g\in G}\lambda^2_{pq}(g)=1$, where $\lambda_{pq}(g)\geq0$,  it follows
that there exists exactly one $h_{pq}\in G$
such that
\[
\lambda_{pq}(h_{pq})=1,
\qquad
\lambda_{pq}(g)=0
\quad\text{if }g\ne h_{pq}.
\]
Since the above discussion holds for arbitrary
$p,q\in[n]$, we derive
\[
(Z_\pi)_{pq}=\pi(h_{pq}),
\qquad \forall\pi\in\widehat G,
\]
and $\hat Z_g\in\{0,1\}^{n\times n}$ for every $g\in G$.

Next, we are going to show there exist $g = (g_1,\dots,g_n)\in G^n$ such that 
\[
Z_\pi = V_\pi(g) V_\pi(g)^H, \qquad\forall\pi\in\widehat G\setminus\{\pi_0\}.
\]
Since $Z_\pi=V_\pi V_\pi^H$
and each $V_{\pi,i}$ is unitary, for any
$p,q,t\in[n]$, we have
\[
\begin{aligned}
\pi(h_{pq}h_{qt})=\pi(h_{pq})\pi(h_{qt}) = (Z_\pi)_{pq}(Z_\pi)_{qt}
&=
V_{\pi,p}V_{\pi,q}^H
V_{\pi,q}V_{\pi,t}^H\\
&=
V_{\pi,p}V_{\pi,t}^H
=(Z_\pi)_{pt}=\pi(h_{pt}).
\end{aligned}
\]
where the first equality comes from the homomorphism property of group representation. Therefore, we derive
\[
\pi(h_{pq}h_{qt})=\pi(h_{pt}),
\qquad \forall\pi\in\widehat G.
\]
Since the irreducible representations of $G$
separate group elements\footnote{See \Cref{append:irreducible}}, we obtain
\[
h_{pq}h_{qt}=h_{pt}.
\]
Fixing $t=1$ and using $h_{1q}=h_{q1}^{-1}$,
we derive
\[
h_{pq}=h_{p1}h_{q1}^{-1}.
\]
Let $g_p:=h_{p1}$ for every $p\in[n]$.
Then
\[
h_{pq}=g_pg_q^{-1},
\qquad \forall p,q\in[n].
\]
Thus, for every $\pi\in\widehat G\setminus\{\pi_0\}$,
we have
\[
(Z_\pi)_{pq}
=
\pi(g_pg_q^{-1})
=
\pi(g_p)\pi(g_q)^H.
\]
This implies that
\[
Z_\pi=V_\pi(g)V_\pi(g)^H,
\qquad
\forall\pi\in\widehat G\setminus\{\pi_0\},
\]
where $g=(g_1,\ldots,g_n)\in G^n$.
Hence $\Z\in\cQ_G^n$.

($\Rightarrow$) Conversely, if $\Z\in\cQ_G^n$, then
$Z_\pi=V_\pi(g)V_\pi(g)^H$ for some $g\in G^n$.
Since
$V_\pi(g)^HV_\pi(g)=nI_{d_\pi}$,
we have $\rank(Z_\pi)=d_\pi$ for every
$\pi\in\widehat G\setminus\{\pi_0\}$.
\end{proof}

\begin{rmk}
If $G$ is a finite abelian group, then $d_\pi=1$ for every $\pi\in\widehat G\setminus\{\pi_0\}$. In this case, $\Z$ has minimal rank if and only if each $Z_\pi$ has rank one.
\end{rmk}

Before proceeding to discuss our main result, we start with the
following quantities. Suppose $g=(g_1,\dots,g_n)\in G^n$.
We consider a map
$\varphi:(\RR^{n\times n})^{|G|}\times G^n
\rightarrow(\RR^{n\times n})^{|G|}$,
$(\boldsymbol{F},g)\mapsto
\boldsymbol{M}=(M_h)_{h\in G}$,
where
\[
(M_h)_{ij}
=
(F_{g_i h g_j^{-1}})_{ij}
-
(F_{g_i g_j^{-1}})_{ij},
\qquad \forall h\in G.
\]
Then for each $h\in G$, denote
\begin{equation}\label{eq:u_k}
u_h=M_h\vone_n,
\end{equation}
and for each $\pi\in\widehat G$, denote
\begin{equation}
\label{eq:zkstar}
V_\pi(g):=
\begin{pmatrix}
\pi(g_1)\\
\vdots\\
\pi(g_n)
\end{pmatrix},
\qquad
\Pi_\pi:=I_{nd_\pi}
-\frac1n V_\pi(g)V_\pi(g)^H.
\end{equation}
For simplicity, denote $\Pi:=\Pi_{\pi_0} = I_n - \tfrac{1}{n}\vone_n \vone_n^T$, where
$\pi_0$ is the trivial representation. Then we define
\begin{equation}\label{def:alpha_beta}
\begin{aligned}
\alpha_g
&:=\underset{h\neq e,\;j\in[n]}{\min}(u_h)_j,
\quad
\beta_g
:=\underset{h\neq e}{\max}
\tfrac{1}{\sqrt n}\|\Pi u_h\|_2.
\end{aligned}
\end{equation}
and
\begin{equation}\label{def:rho}
\rho_g:=
\underset{\pi\in\widehat G\setminus\{\pi_0\}}{\max}
\big(-\lambda_{\min}(\Pi_\pi C_\pi\Pi_\pi)\big)_+.
\end{equation}

Based on the above quantities, we are ready to show the main
exactness result for SDP relaxation \eqref{eq:sdp}.

\begin{thm}\label{thm:exactness}
    Given $g\in G^n$, if we have
    \begin{equation}\label{eq:exact_condition}\tag{$*$}
        \alpha_g > \rho_g + 2\sqrt{|G|-1}\,\beta_g,
    \end{equation}
    then the SDP \eqref{eq:sdp} is tight, and gives the
    optimal solution corresponding to $g$.
\end{thm}

\subsection{Dual multipliers}

Before proceeding to the proof of the above theorem, we introduce the dual problem of \eqref{eq:sdp}.  Denote by \(\prod_{\pi\in\widehat G}\mathbb H^{nd_\pi}\) the Cartesian product of Hermitian matrix spaces, where each element is a tuple \(\Y=(Y_\pi)_{\pi\in\widehat G}\) with \(Y_\pi\in\mathbb H^{nd_\pi}\) for every \(\pi\in\widehat G\). The dual problem of \eqref{eq:sdp} is defined as follows,
\begin{equation}\label{eq:sdp-dual}\tag{D}
\begin{split}
     \max_{\Y=(Y_\pi)_{\pi\in\widehat G}}\quad&
     \sum_{\pi\in\widehat G}d_\pi\tr(Y_\pi)
     - E_n\bullet Y_{\pi_0}\\
     \text{s.t.}\quad&
     C_\pi-Y_\pi\succeq\vzero,\quad
     \forall\pi\in\widehat G\setminus\{\pi_0\},\\
     &\Y\in\cD_G^n
\end{split}
\end{equation}
where
\begin{equation}\label{eq:D_L}
\begin{split}
\cD_G^n:=\Big\{
\Y=(Y_\pi)_{\pi\in\widehat G}&\in
\prod_{\pi\in\widehat G}\mathbb H^{nd_\pi}:\\
&\hat{\Y}=(\hat Y_h)_{h\in G}
=\cF_G^{-1}((Y_\pi)_{\pi\in\widehat G}),\\
&\offdiag(\hat Y_h)\geq 0,\quad \forall h\in G,\\
&\hat Y_h\in\RR^{n\times n},\quad \forall h\in G
\Big\}.
\end{split}
\end{equation}

Then we introduce some definitions and notation that are useful in the proof. If \eqref{eq:exact_condition} holds, there exists $\varepsilon>0$ such that
\begin{equation}\label{eq:b}
b:=\alpha_g-\sqrt{|G|-1}\beta_g-\varepsilon
>\rho_g+\sqrt{|G|-1}\beta_g.
\end{equation}
For any $h\neq e$, let $v_h:=u_h-b\vone_n$, where $u_h$ is defined in \eqref{eq:u_k}, and $\tau_h=\vone_n^Tv_h$. Since $M_{h^{-1}}=M_h^T$, we have $\tau_h=\tau_{h^{-1}}$. Let
\[
N_h:=\frac{1}{\tau_h}v_hv_{h^{-1}}^T.
\]
Then it is true that
\begin{equation}\label{eq:tau}
v_h\geq(\sqrt{|G|-1}\beta_g+\varepsilon)\vone_n>\vzero,
\qquad
\tau_h\geq n(\sqrt{|G|-1}\beta_g+\varepsilon).
\end{equation}
It is straightforward that $N_h>\vzero$ for every $h\neq e$. Moreover,
\[
N_{h^{-1}}=N_h^T,\qquad N_h\vone_n=v_h.
\]
Denote $N_e=\vzero$. Then we construct the dual multipliers in \eqref{eq:D_L}, $\hat{\Y}:=(\hat Y_h)_{h\in G}$, in the following way,
\begin{equation}\label{dual_multiplier}
(\hat Y_h)_{ij}
=
(N_{g_i^{-1}hg_j})_{ij}
-b\,\mathds{1}_{\{h=e\}}\mathds{1}_{\{i=j\}},
\qquad \forall i,j\in[n],\ h\in G.
\end{equation}
It follows that $\hat Y_{h^{-1}}=\hat Y_h^T$ and
$\offdiag(\hat Y_h)\geq0$ for every $h\in G$.

Then we discuss some properties of the dual multipliers constructed as above.

\begin{lem}\label{lem:prop_dual_1}
Consider $\hat{\Y}$ as constructed in \eqref{dual_multiplier},
and let $\Y=\mathcal F_G(\hat{\Y})$. Then for every
$\pi\in\widehat G$, we have
\begin{enumerate}[label=(\roman*)]
\item $(Y_\pi V_\pi(g))_i
=\pi(g_i)\big(\sum_{h\neq e}(v_h)_i\pi(h)-bI_{d_\pi}\big)$,
$\forall i\in[n]$;
\item $(Y_\pi)_{ij}
=\sum_{h\neq e}\frac{(v_h)_i(v_{h^{-1}})_j}{\tau_h}
\pi(g_ihg_j^{-1})
-b\mathds{1}_{\{i=j\}}I_{d_\pi}$,
$\forall i,j\in[n]$.
\end{enumerate}
\end{lem}

\begin{proof}
$(i)$: From \eqref{dual_multiplier}, for any $i,j\in[n]$,
we have
\begin{equation}\label{eq:prop_dual}
\begin{split}
(Y_\pi)_{ij}
&=\sum_{t\in G}(\hat Y_t)_{ij}\pi(t)\\
&=\sum_{t\in G} (N_{g_i^{-1}tg_j})_{ij}\pi(t)
-b \mathds{1}_{\{i=j\}}I_{d_\pi}\\
&=\sum_{h\in G}(N_h)_{ij}
\pi(g_ihg_j^{-1})
-b\mathds{1}_{\{i=j\}}I_{d_\pi}\\
&=\pi(g_i)\Big(\sum_{h\neq e}
(N_h)_{ij}\pi(h)\Big)\pi(g_j)^H
-b\mathds{1}_{\{i=j\}}I_{d_\pi}.
\end{split}
\end{equation}
which implies
\[
\begin{aligned}
(Y_\pi V_\pi(g))_i
&=\sum_{j=1}^n(Y_\pi)_{ij}\pi(g_j)\\
&=\pi(g_i)\sum_{h\neq e}\pi(h)
\Big(\sum_{j=1}^n(N_h)_{ij}\Big)
-b\pi(g_i)\\
&=\pi(g_i)\sum_{h\neq e}
\Big((N_h\vone_n)_i\Big)\pi(h)-b\pi(g_i)\\
&=\pi(g_i)\Big(
\sum_{h\neq e}(v_h)_i\pi(h)-bI_{d_\pi}
\Big).
\end{aligned}
\]
Thus, $(i)$ is proved. Then we show $(ii)$ is true.
From \eqref{eq:prop_dual}, we have
\[
\begin{aligned}
(Y_\pi)_{ij}
&=\sum_{h\neq e}(N_h)_{ij}
\pi(g_ihg_j^{-1})
-b\mathds{1}_{\{i=j\}}I_{d_\pi}\\
&=\sum_{h\neq e}
\frac{(v_h)_i(v_{h^{-1}})_j}{\tau_h}
\pi(g_ihg_j^{-1})
-b\mathds{1}_{\{i=j\}}I_{d_\pi},
\end{aligned}
\]
which proves $(ii)$.
\end{proof}

\subsection{Proof of \Cref{thm:exactness}}
Now we are ready to prove the above theorem.

\begin{proof}
Suppose that \eqref{eq:exact_condition} holds. There exist $b$ and $\hat{\Y}:=(\hat Y_h)_{h\in G}$ as constructed in \eqref{eq:b} and \eqref{dual_multiplier}. Let $\Y = \mathcal{F}_G(\hat{\Y})$. We consider the following dual slack matrices,
\[
Q_\pi := C_\pi - Y_\pi,\quad \pi\in\widehat G\setminus\{\pi_0\}.
\]
It is enough to show that for any $\pi\in\widehat G\setminus\{\pi_0\}$,
\begin{equation}
    Q_\pi\succeq \vzero,\qquad\rank(Q_\pi)=(n-1)d_\pi,\qquad Q_\pi V_\pi(g) = \vzero.
\end{equation}
We first show that $Q_\pi V_\pi(g)=0$. From the definition of $u_h$, it is true that
\[
\begin{aligned}
(C_\pi V_\pi(g))_i
&=\sum_{j=1}^n(C_\pi)_{ij} \pi(g_j)\\
&=\sum_{j=1}^n \Big(\sum_{t\in G} (F_t)_{ij} \pi(t)\Big) \pi(g_j)\\
&=\sum_{j=1}^n \sum_{t\in G} (F_t)_{ij} \pi(tg_j)
\end{aligned}
\]
By taking $t = g_i h g_j^{-1}$, we derive
\[
\begin{aligned}
(C_\pi V_\pi(g))_i &= \sum_{j=1}^n \sum_{h\in G} (F_{g_i h g_j^{-1}})_{ij} \pi(g_i) \pi(h)\\
&= \pi(g_i)\Big(\sum_{h\in G}(u_h)_i\pi(h)
+\sum_{j=1}^n(F_{g_i g_j^{-1}})_{ij}\sum_{h\in G}\pi(h)\Big)\\
&= \pi(g_i)\sum_{h\neq e}(u_h)_i\pi(h).
\end{aligned}
\]
where the second equality comes from the definition that $(F_{g_i h g_j^{-1}})_{ij} = (u_h)_i + (F_{g_i g_j^{-1}})_{ij}$ as in \eqref{eq:u_k}, and the last equality is from the fact that $\sum_{h\in G}\pi(h)=0$ for every nontrivial irreducible representation $\pi$.\footnote{See \Cref{append:orthogonality}. When $G = \ZZ_L$, we have $\pi_k(h)=\omega^{-kh},\forall k,h\in \ZZ_L$. In this case, $\sum_{h\in G}\pi(h)=0$ is equivalent to $\sum_{h=0}^{L-1}\omega^{-kh}=0, k \in \ZZ_L \setminus \{0\}$.}
From \Cref{lem:prop_dual_1}, for every $i\in[n]$, we have
\begin{equation*}
\begin{split}
(Y_\pi V_\pi(g))_i &= \pi(g_i)\Big(\sum_{h\neq e}(v_h)_i\pi(h)-bI_{d_\pi}\Big).
\end{split}
\end{equation*}
It follows that
\begin{equation*}
\begin{split}
(Q_\pi V_\pi(g))_i
&= \pi(g_i)\Big(bI_{d_\pi}+\sum_{h\neq e}((u_h)_i-(v_h)_i)\pi(h)\Big)\\
&= b\pi(g_i)\Big(I_{d_\pi}+\sum_{h\neq e}\pi(h)\Big)\\
&= b\pi(g_i)\sum_{h\in G}\pi(h)=0.
\end{split}
\end{equation*}
where the second equality comes from how $v_h$ is constructed, that is, $v_h = u_h - b \vone_n$, and the last equality is from the fact that $\sum_{h\in G}\pi(h)=0$ for every nontrivial irreducible representation $\pi$.

Next, we show $Q_\pi\succeq\vzero$ and $\rank(Q_\pi)=(n-1)d_\pi$. From \Cref{lem:prop_dual_1}, we have
\[
Y_\pi+bI_{nd_\pi}
=\sum_{h\neq e}\frac{1}{\tau_h}W_{\pi,h}\pi(h)W_{\pi,h^{-1}}^H,
\]
where
\[
W_{\pi,h}:=
\begin{pmatrix}
(v_h)_1\pi(g_1)\\
\vdots\\
(v_h)_n\pi(g_n)
\end{pmatrix}.
\]
Since we have
\[
\begin{aligned}
\Pi_\pi W_{\pi,h}
&=\left(I_{nd_\pi}-\frac{1}{n}V_\pi(g)V_\pi(g)^H\right)W_{\pi,h}\\
&=W_{\pi,h}-\frac{1}{n}V_\pi(g)(\vone_n^Tv_h)\\
&=\begin{pmatrix}
\big((v_h)_1 - \tfrac{1}{n}\vone_n^T v_h\big)\pi(g_1)\\
\vdots\\
\big((v_h)_n - \tfrac{1}{n}\vone_n^T v_h\big)\pi(g_n)
\end{pmatrix}\\
&=
\begin{pmatrix}
(\Pi v_h)_1\pi(g_1)\\
\vdots\\
(\Pi v_h)_n\pi(g_n)
\end{pmatrix}
\end{aligned}
\]
it follows that
\begin{equation}\label{eq:rho}
    \|\Pi_\pi W_{\pi,h}\|
    =\|\Pi u_h\|_2
    \leq \sqrt{n}\beta_g.
\end{equation}
Combining \eqref{eq:tau} and \eqref{eq:rho}, it follows that
\[
\begin{aligned}
\|\Pi_\pi(Y_\pi+bI_{nd_\pi})\Pi_\pi\|
&\leq \sum_{h\neq e}
\frac{\|\Pi_\pi W_{\pi,h}\|
\|\Pi_\pi W_{\pi,h^{-1}}\|}{\tau_h}\\
&\leq \sum_{h\neq e}\frac{n\beta_g^2}{\tau_h}\\
&\leq \sum_{h\neq e}\frac{\beta_g^2}
{\sqrt{|G|-1}\beta_g+\varepsilon}\\
&\leq \frac{(|G|-1)\beta_g^2}
{\sqrt{|G|-1}\beta_g+\varepsilon}
\leq \sqrt{|G|-1}\beta_g.
\end{aligned}
\]
where the second last inequality comes from \eqref{eq:tau}.
From the definition of $\rho_g$, we have that $\forall w\in\langle V_\pi(g)\rangle^\perp \setminus \{0\}$,
\[
w^H C_\pi w\geq-\rho_g\|w\|_2^2.
\]
Then it follows that
\[
\begin{aligned}
w^H Q_\pi w
&=w^H(C_\pi-Y_\pi)w\\
&=w^H\big(C_\pi+bI_{nd_\pi}-(Y_\pi+bI_{nd_\pi})\big)w\\
&\geq(b-\rho_g-\sqrt{|G|-1}\beta_g)\|w\|_2^2>0.
\end{aligned}
\]
which implies that $\rank(Q_\pi) \geq (n-1) d_\pi$. Combining the above results with $Q_\pi V_\pi (g) = \vzero$, it is true that $Q_\pi\succeq\vzero$ and $\rank(Q_\pi)=(n-1)d_\pi$, $\forall\pi\in\widehat G\setminus\{\pi_0\}$. Thus, the proof is finished.
\end{proof}

\section{Multireference Alignment}\label{sec:mra}

In this section, we apply \eqref{eq:sdp} to the multireference
alignment (MRA) problem. Using the sufficient condition established earlier, we analyze the exactness of the SDP relaxation and show
that it achieves the minimax-optimal noise threshold for exact recovery.

\subsection{Multireference alignment}
Let $G\cong\mathbb{Z}_L$ be a cyclic group acting on
$\RR^{m\times L}$, and let $g_0$ denote its generating element,
which acts by cyclically shifting the columns of a matrix by one
position. We use the convention $(g_0^s x)(\ell)=x(\ell+s)$,
where $x(\ell)$ is the $\ell$-th column of $x$ and the column indices are taken modulo $L$.
Thus, each $g\in\mathcal{G}$ acts by cyclically shifting the columns by some number of positions. Consider the \textit{multireference alignment} (MRA) problem of estimating an unknown signal $x\in \RR^{m\times L}$ from measurements of the form:
\begin{equation*}
    y_i =  g_i x + \sigma\epsilon_i,\quad i =1,2\dots,n
\end{equation*}
where $ g_1,\dots, g_n$ are unknown elements from $\mathcal{G}\cong \mathbb{Z}_L$, and the noise entries satisfy
$(\epsilon_i)_{k\ell}
\overset{\mathrm{i.i.d.}}{\sim}\mathcal{N}(0,1)$. To recover the group actions, we consider the following minimization problem
\begin{equation}\label{eq:mra}
    \min_{s\in \ZZ^n_L} \sum_{i,j\in [n]} \|g_0^{-s_i}y_i - g_0^{-s_j}y_j\|^2_F 
\end{equation}
The underlying signal \(x\) can be estimated by taking average of those aligned measurements. Equivalently, define the potential function
\[
f_{ij}(t) = \frac{1}{2}\|g_0^{-t}y_i - y_j\|_F^2,
\]
then \eqref{eq:mra} becomes
\[
\min_{s\in\ZZ^n_L} \sum_{i,j\in [n]} f_{ij}(s_i-s_j)
\]
Thus, multireference alignment is an instance of \eqref{eq:cyclic-difference}.

Suppose that $g_i = g_0^{s_i^*}$,that is, $y_i = g_0^{s_i^*} x + \sigma \epsilon_i, i \in [n]$
where \(s^*\in\mathbb Z_L^n\) is the unknown shift vector. Let \(d_{ij}:=s^*_i-s^*_j, \forall 1\leq i,j\leq n\). For \(t\in\ZZ_L\), define 
\((F_t)_{ij}:=f_{ij}(t)=\tfrac{1}{2}\|y_i\|_F^2 + \tfrac{1}{2}\|y_j\|_F^2 - \langle g_0^{-t}y_i,y_j\rangle\). The objective matrix in \eqref{eq:multi-frequency-nonconstant} is 
\begin{equation*}
\C := (C_0,\dots,C_{L-1}) = \cF_{\ZZ_L}(F_0,\dots,F_{L-1}) 
\end{equation*}

Let $\delta_q:=\|g_0^{-q}x-x\|_F, q\in\ZZ_L$,  $\delta_{\min} = \min_{q\neq 0} \delta_q$, and $\delta_{\max}:=\max_{q\neq 0} \delta_q$.  Let \(e_i:=g_0^{-s_i^*}\epsilon_i,~\forall i\in[n]\), and $\bar e := \tfrac{1}{n} \sum_{i=1}^n e_i$. Denote $\tilde y_i := g_0^{-s_i^*}y_i = x + \sigma e_i,\;\forall i\in [n]$.

\subsection{Exactness results}

We provide a noise threshold for exact recovery in the MRA model \eqref{eq:mra}.

\begin{thm}\label{thm:mra_minimax}
Suppose that $L\geq 2$ is fixed and $m = o(n \log n)$. If $\delta_{\max}/\delta_{\min}$ is bounded by a constant and the noise level $\sigma$ satisfies
\[
\sigma \leq (1-\eta) \frac{\delta_{\min}}{2\sqrt{2 \log n}}
\]
where $\eta \in (0,1)$ is a constant, then \eqref{eq:sdp} is tight w.h.p.
\end{thm}

It is enough to prove that under the above assumptions \eqref{eq:exact_condition} holds with high probability. Throughout this section, for simplicity, we denote $s := s^*$. We discuss the orders of $\alpha_s$, $\beta_s$ and $\rho_s$ which are defined in \eqref{def:alpha_beta} and \eqref{def:rho}. We first show that $\alpha_s/n$ has a constant lower bound with high probability.

\begin{lem}\label{lem:mra_alpha}
    There exists a constant $c_\eta > 0$ 
    such that $\alpha_s/n \geq c_\eta \delta^2_{\min}$ w.h.p.
\end{lem}
\begin{proof}
Denote $\bar y:= x + \sigma \bar e$ and $a_q:=(I-g_0^{q})\bar y$. Since it is true that
\[
(F_{d_{ij}+q})_{ij} - (F_{d_{ij}})_{ij} = \langle \tilde y_i, \tilde y_j - g_0^{q} \tilde y_j\rangle
\]
it follows that 
\begin{equation}\label{eq:mra_alpha_key}
    \frac{(u_q)_i}{n} = \langle \tilde y_i, a_q \rangle = \langle \bar y, a_q \rangle + \sigma \langle e_i - \bar e, a_q \rangle = \frac{1}{2}\|a_q\|_F^2 +  \sigma \langle e_i - \bar e, a_q \rangle.
\end{equation}
where the last equality comes from $\langle \bar y, a_q \rangle = \langle \bar y, \bar y - g_0^q \bar y \rangle = \tfrac{1}{2}\| \bar y - g_0^q \bar y\|_F^2$.

Then we discuss the two terms $\|a_q\|_F^2$ and $ \langle e_i - \bar e, a_q \rangle$. Since we have $\bar e \sim \tfrac{1}{\sqrt{n}}\mathcal{N}(0, I_{mL})$, it is true that $\|\bar e\|_F^2 = O_{\mathbb{P}}(\tfrac{mL}{n})$. Hence we obtain
\begin{equation}\label{eq:mra_alpha_first_term}
    \|a_q\|_F = \delta_{q} + O_{\mathbb P}\Big(\sigma \sqrt{\tfrac{mL}{n}}\Big) = \delta_{q} + o_{\mathbb P}(\delta_{q})
\end{equation}
where the last equality comes from the assumption that $L$ is fixed and $\tfrac{m}{n \log n}=o(1)$. It remains to show that for the given constant $\eta\in(0,1)$, there exists a constant $c_\eta\in(0,1)$ such that with high probability,
\begin{equation}\label{eq:mra_alpha_second_term}
    \sigma \langle e_i - \bar e, a_q \rangle \geq - \tfrac{1}{2}(1 - c_\eta) \|a_q\|_F^2, \;\forall\; i \in [n]
\end{equation}
Since $e_i-\bar e$ and $\bar e$ are independent, we get 
\[
e_i-\bar e \;|\; \bar e \sim \mathcal{N}(0,(1-\tfrac{1}{n})I_{mL})
\]
Since $a_q$ depends only on $\bar e$, it follows that
\[
\langle e_i - \bar e, a_q/\|a_q\|_F \rangle \sim \mathcal{N}(0,1-1/n)
\] 
Thus,  for any fixed $\epsilon > 0$, with probability at least $1-n^{-(1+\epsilon)}$,
\[
\langle e_i - \bar e, a_q/\|a_q\|_F \rangle \geq -\sqrt{2(1+\epsilon)\log n} 
\]
Consequently, using the union bound, with probability at least $1-(L-1)n^{-\epsilon}$, 
\[
\langle e_i - \bar e, a_q/\|a_q\|_F \rangle \geq -\sqrt{2(1+\epsilon)\log n}, \; \forall\;  i \in [n].
\]
Given $\eta\in(0,1)$, we can choose $\epsilon>0$ and
$c_\eta\in(0,1)$ such that
\[
(1-\eta)\sqrt{1+\epsilon}<1-c_\eta.
\]
and
\eqref{eq:mra_alpha_second_term} holds. Combining \eqref{eq:mra_alpha_key}, \eqref{eq:mra_alpha_first_term} and \eqref{eq:mra_alpha_second_term}, we derive
\[
\begin{aligned}
\frac{(u_q)_i}{n}
&\geq \frac12\|a_q\|_F^2
-\frac12(1-c_\eta)\|a_q\|_F^2\\
&=\frac{c_\eta}{2}\|a_q\|_F^2\\
&\geq\frac{c_\eta}{2}(1-o_{\mathbb P}(1))\delta_{\min}^2.
\end{aligned}
\]
Thus, the proof is finished. 
\end{proof}

Then we show that $\beta_s/(n\delta^2_{\min})$ vanishes as $n \rightarrow \infty$. 
\begin{lem}\label{lem:mra_beta}
    It is true that $\beta_s/n = o_{\mathbb P}(\delta_{\min}^2)$.
\end{lem}
\begin{proof}
From \eqref{eq:mra_alpha_key}, we have
\[
u_q = \frac{n}{2} \|a_q\|_F^2\vone_n + \sigma n \;\big( \vct(e_i - \bar e)^T \vct(a_q)  \big)
\]
Consider a matrix which is defined as follows
\[
G : =\begin{pmatrix}
    \vct(e_1 - \bar e) & \cdots & \vct(e_n - \bar e)
\end{pmatrix} = E\Pi
\]
where $E$ is a matrix which has independent standard Gaussian columns. It follows that
\[
\|G\| \leq \|E\| = O_{\mathbb P}(\sqrt{n} + \sqrt{mL})
\]
and
\[
\Pi u_q = \Pi \Big(\frac{n}{2} \|a_q\|_F^2\vone_n + \sigma n \;\big( \vct(e_i - \bar e)^T \vct(a_q)  \big)\Big) = n\sigma \; G^T \vct(a_q)
\]
Since we have $\|G\| = O_{\mathbb{P}}(\sqrt{n} + \sqrt{mL})$, and $\|a_q\|_F = \delta_q + o_{\mathbb P}(\delta_{\min})$, it follows that
\[
\frac{\beta_s}{n\delta^2_{\min}} \leq \frac{1}{\delta^2_{\min}}\frac{\sigma}{\sqrt{n}} \|G\| \big(\max_{q\neq 0}  \|a_q\|_F \big) = O_{\mathbb P}\Big(\frac{1}{\sqrt{\log n}}+\sqrt{\frac{m}{n \log n}}\Big) = o_{\mathbb P}(1)
\]
where the last equality comes from assumption that there exists a constant $K$ such that $\delta_{\max}/\delta_{\min} \leq K$.
\end{proof}

Finally, we prove that $\rho_s/(n\delta^2_{\min})$ also vanishes as $n\rightarrow \infty$. 
\begin{lem}\label{lem:mra_rho}
    It is true that $\rho_s/n = o_{\mathbb P}(\delta_{\min}^2)$.
\end{lem}
\begin{proof}
Denote 
\begin{equation*}
 z^{\odot k} :=
 \begin{pmatrix}
 \omega^{-ks_1}\\
 \vdots\\
 \omega^{-ks_n}
 \end{pmatrix},
 \qquad
 \Pi_k:=I-\frac1n z^{\odot k}(z^{\odot k})^H.
\end{equation*}
which is a special case of \eqref{eq:zkstar} when $G = \ZZ_L$.
For any $i \in [n]$, let 
\[
(\hat{y}_i(1),\dots,\hat{y}_i(L)) = \cF_{\ZZ_L}(y_i(1),\dots,y_i(L))\]
where $y_i(j)$ is the $j$-th column of $y_i$. Let 
$W_k:= \begin{pmatrix}
    \hat y_1(k) & \cdots & \hat y_n(k)
\end{pmatrix}$.
For $k\neq0$, since $\sum_{t=0}^{L-1}\omega^{kt}=0$,
we have
\[
\begin{aligned}
(C_k)_{ij}
&=-\sum_{t=0}^{L-1}\omega^{-kt}
\langle g_0^{-t}y_i,y_j\rangle\\
&=-\sum_{t=0}^{L-1}\sum_{\ell=0}^{L-1}
\omega^{-kt}y_i(\ell-t)^T y_j(\ell)\\
&=-\sum_{r=0}^{L-1}\sum_{\ell=0}^{L-1}
\omega^{-k(\ell-r)}y_i(r)^T y_j(\ell)\\
&=-\hat y_i(k)^H\hat y_j(k).
\end{aligned}
\]
Therefore, we have $C_k=-W_k^HW_k$ for $k\neq0$. 

Denote $(\hat{x}(1),\dots,\hat{x}(L)) := \cF_{\ZZ_L}(x(1),\dots,x(L))$ and 
\[
E_k :=
    (\hat \epsilon_1(k),  \cdots,  \hat \epsilon_n(k))
 = \cF_{\ZZ_L}(\epsilon_1(k) , \cdots, \epsilon_n(k)),\quad \forall k \in [L-1].
 \]
Since $\hat{y}_i(k) = \omega^{ks_i} \hat{x}(k) + \sigma \hat{\epsilon}_i(k)$, we also have
\[
W_k = \hat x(k) (z^{\odot k})^H + \sigma E_k
\]
Then we obtain that 
\[
\Pi_k C_k \Pi_k = -\sigma^2 \Pi_kE_k^HE_k \Pi_k
\]
which implies that
\[
\rho_s \leq \sigma^2 \max_{k\neq 0}\|E_k\Pi_k\|^2
\]
Recall the entries of each $\epsilon_i$ are i.i.d. standard Gaussian, and
\[
\hat\epsilon_i(k)
=\sum_{\ell=0}^{L-1}
\left[
\cos\left(\frac{2\pi k\ell}{L}\right)
-\boldsymbol{\mathrm{i}}\sin\left(\frac{2\pi k\ell}{L}\right)
\right]\epsilon_i(\ell).
\]
If $2k\not\equiv 0\bmod L$, then
\begin{equation}\label{eq:sin_cos}
\sum_{\ell=0}^{L-1}\cos^2\left(\frac{2\pi k\ell}{L}\right)
=\sum_{\ell=0}^{L-1}\sin^2\left(\frac{2\pi k\ell}{L}\right)
=\frac L2,
\end{equation}
and
\begin{equation}\label{eq:independent_sin_cos}
\sum_{\ell=0}^{L-1}
\cos\left(\frac{2\pi k\ell}{L}\right)
\sin\left(\frac{2\pi k\ell}{L}\right)=0.
\end{equation}
Thus, from $\eqref{eq:sin_cos}$ and \eqref{eq:independent_sin_cos} the real and imaginary parts of $E_k$ are
independent Gaussian matrices with i.i.d. $\cN(0,L/2)$ entries,  and
\[
E_k\overset{d}{=}\sqrt{\frac L2}(G_1+\mathrm{i}G_2),
\]
where $G_1,G_2\in\RR^{m\times n}$ are independent standard Gaussian
random matrices. When $L$ is even and $k=L/2$, the Fourier weights
are $(-1)^\ell$, and hence
\[
E_k\overset{d}{=}\sqrt L G_1.
\]
In each case, we have
\[
\|E_k\|=O_{\mathbb P}\left(\sqrt L(\sqrt m+\sqrt n)\right).
\]
Consequently, it follows that
\[
\frac{\rho_s}{n\delta^2_{\min}} \leq  \frac{\sigma^2}{n\delta_{\min}^2} \max_{k\neq 0}\|E_k\|^2 = O_{\mathbb P}\Big(\frac{1}{\log n} + \frac{m}{n \log n}\Big) = o_{\mathbb P}(1) \qedhere
\]
\end{proof}

\begin{proof}[Proof of \Cref{thm:mra_minimax}]
    From \Cref{lem:mra_alpha} -- \ref{lem:mra_rho}, we have
    \[
    \frac{\alpha_s}{n} \geq c_\eta \delta_{\min}^2 > \frac{\rho_s}{n} + \frac{2\sqrt{L-1}\beta_s}{n} = o_{\mathbb P}(\delta_{\min}^2). 
    \]
    Thus, \eqref{eq:exact_condition} holds w.h.p.
\end{proof}

\subsection{Minimax optimality} We conclude this section by establishing a converse result.
Specifically, we show that if the noise level exceeds the
recovery threshold established above, then no estimator can recover the true shifts, up to
a global cyclic shift, with high probability.

\begin{thm}
\label{thm:mra-converse}
For every fixed $\eta>0$, if the noise level $\sigma$ satisfies 
\begin{equation*}
\label{eq:mra-threshold-above}
\sigma
\ge
(1+\eta)
\frac{\delta_{\min}}{2\sqrt{2\log n}},
\end{equation*}
then no method can recover the true solution up to a global shift.
\end{thm}

\begin{proof}
See \Cref{append_proof_of_mra_converse}.
\end{proof}

\section{Joint Assignment}\label{sec:joint_assignment}

In this section, we consider another application: joint assignment problem. We apply our exactness condition to this problem and show that \eqref{eq:sdp} achieves minimax optimality.

\subsection{Joint Assignment}
There are unknown planted labels
\[
s_i^*\in\mathbb Z_L,
\qquad i=1,\ldots,n,
\]
and for some observed pairs \((i,j)\in\Omega\), we observe noisy relative measurements
\[
Y_{ij}\approx s_i^*-s_j^*\bmod L.
\]

We assume that each unordered pair $(i,j)$ is observed independently with probability $p$, and the measurements are independent across observed pairs. For each observed pair, we set $Y_{ji}=-Y_{ij}\pmod L$.

A likelihood-based estimator chooses labels \(s\in\mathbb Z_L^n\) that best explain the pairwise measurements:
\begin{equation}
\max_{s\in\mathbb Z_L^n}
\sum_{(i,j)\in\Omega}
\log p_{ij}\bigl(Y_{ij}\mid s_i-s_j\bigr).
\label{eq:joint-assignment-mle}
\end{equation}
Equivalently, defining the negative log-likelihood potential
\[
f_{ij}(t)
=
-\log p_{ij}(Y_{ij}\mid t),
\qquad t\in\mathbb Z_L,
\]
the maximum likelihood estimator becomes
\begin{equation}
\min_{s\in\mathbb Z_L^n}
\sum_{(i,j)\in\Omega}
f_{ij}(s_i-s_j).
\label{eq:joint-assignment-min}
\end{equation}
Thus joint assignment is also an instance of the general cyclic-difference problem \eqref{eq:cyclic-difference}.

In a random-corruption model, one observes
\[
Y_{ij}=
\begin{cases}
s_i^*-s_j^* \pmod L, & \text{with probability }\pi_0,\\
\text{a uniformly random element of }\mathbb Z_L, & \text{with probability }1-\pi_0.
\end{cases}
\]
Let $q_0=\pi_0+\tfrac{1-\pi_0}{L}$, and $q_1=\tfrac{1-\pi_0}{L}$. Then
\[
p_{ij}(Y_{ij}\mid t)
=
\begin{cases}
q_0, & t=Y_{ij},\\[2mm]
q_1, & t\neq Y_{ij},
\end{cases}
\]
and hence
\[
f_{ij}(t)
=
-\log p_{ij}(Y_{ij}\mid t).
\]
Equivalently, we define
\[
(F_t)_{ij} = 
\begin{cases}
- \mathds{1}\{(i,j)\in\Omega\mid Y_{ij} = t\},\quad &\text{if}\;\; i \neq j\\
0,\quad &\text{if}\;\; i = j
\end{cases}
\]
The objective matrix in \eqref{eq:multi-frequency-nonconstant} is 
\begin{equation*}
\C := (C_0,\dots,C_{L-1}) = \cF_{\ZZ_L}(F_0,\dots,F_{L-1}). 
\end{equation*}

\subsection{Exactness results}
We provide a threshold on $\pi_0$ under the
random corruption model that guarantees the exactness of
the SDP relaxation.

\begin{thm}\label{thm:exact_joint}
Fix $L\geq 2$ and $\varepsilon>0$. Suppose that $p\geq c_{L,\varepsilon}\frac{\log n}{n}$ for some sufficiently large constant $c_{L,\varepsilon}>0$. If $\pi_0$ satisfies
\begin{equation}\label{eq:joint_pi}
\pi_0 \geq 2\sqrt{\frac{(1+\varepsilon) \log n}{Lpn}}
\end{equation}
then SDP \eqref{eq:sdp} is tight with high probability.
\end{thm}

It is enough to show that if $\pi_0$ satisfies \eqref{eq:joint_pi}, then \eqref{eq:exact_condition} holds with high probability. Throughout this section, for simplicity, we denote $s := s^*$. Let $X_{ij}^{(k)} = (F_{s_i - s_j + k})_{ij} - (F_{s_i - s_j})_{ij}$. For any distinct $i,j\in [n]$, $k\in [L-1]$, we have
\[
X_{ij}^{(k)} = \begin{cases}
    1, \quad &\text{with probability}\quad pq_0\\
    0, \quad &\text{with probability}\quad 1-p(q_0 + q_1)\\
    -1, \quad &\text{with probability}\quad pq_1
\end{cases}
\]
Note that for each fixed \(k\in[L-1]\), the random variables \(\{X_{ij}^{(k)}:1\le i<j\le n\}\) are i.i.d. Then it follows that
\[
\begin{aligned}
(u_k)_i &= \sum_{j\in [n]} (F_{s_i - s_j + k})_{ij} - (F_{s_i - s_j})_{ij}\\
&= \sum_{j \neq i} X_{ij}^{(k)}
\end{aligned}
\]
where the last equality comes from $(F_{t})_{ii}=0, \forall\; i \in [n]$. We first show that $\tfrac{\alpha_s}{np\pi_0}$ has a constant lower bound with high probability.

\begin{lem}\label{lem:joint_alpha}
If $\pi_0$ satisfies \eqref{eq:joint_pi}, there exists a constant $c_\varepsilon>0$ such that
\[
\alpha_s/n > c_\varepsilon p \pi_0 \quad\text{w.h.p.}
\]
\end{lem}

\begin{proof}
From the Chernoff bound, we have 
\[
\begin{aligned}
\mathbb{P}\big((u_k)_i \leq a\big) =\;&  \mathbb{P}\big(e^{-t(u_k)_i} \geq\; e^{-ta}\big)\\
\leq\; & e^{ta} \mathbb{E}\big(e^{-t(u_k)_i}\big)\\
 = \;& e^{ta} \prod_{j\neq i}\mathbb{E}\big(e^{-tX_{ij}^{(k)}}\big)\\
 =\;& e^{ta}\Big(\mathbb{E}\big(e^{-tX_{ij}^{(k)}}\big)\Big)^{n-1}
\end{aligned}
\]
Let $\bar \pi_0 = 2\sqrt{\frac{(1+\varepsilon)\log n}{Lp n}} \leq 2\sqrt{\frac{(1+\varepsilon)}{cL}}$, and denote
\[
\bar q_0=\frac{1+(L-1)\bar\pi_0}{L},
\qquad
\bar q_1=\frac{1-\bar\pi_0}{L}.
\]
Since $\bar q_0>\bar q_1$, by taking $t = \frac{1}{2} \log \frac{\bar q_0}{\bar q_1}>0$, we obtain
\[
\begin{aligned}
\mathbb{E}\big(e^{-tX_{ij}^{(k)}}\big)
&=\; pq_1e^t + pq_0e^{-t} + (1-p(q_0+q_1))\\
&=\; p\bar q_1e^t+p\bar q_0e^{-t}
+(1-p(\bar q_0+\bar q_1))
-p(\pi_0-\bar\pi_0)K\\
&=\; 2p\sqrt{\bar q_0\bar q_1}
+(1-p(\bar q_0+\bar q_1))
-p(\pi_0-\bar\pi_0)K\\
&=\; 1-p(\sqrt{\bar q_0}-\sqrt{\bar q_1})^2
-p(\pi_0-\bar\pi_0)K
\end{aligned}
\]
where
\[
K=\frac{(L-1)(1-e^{-t})+(e^t-1)}{L}.
\]
Thus, it follows that
\[
\begin{aligned}
\mathbb{P}\big((u_k)_i \leq a\big)
&\leq \exp\big\{\tfrac{a}{2}\log \tfrac{\bar q_0}{\bar q_1}\big\}
\Big(1-p(\sqrt{\bar q_0}-\sqrt{\bar q_1})^2
-p(\pi_0-\bar\pi_0)K\Big)^{n-1}\\
&\leq \exp\big\{\tfrac{a}{2}\log \tfrac{\bar q_0}{\bar q_1}
-(n-1)p(\sqrt{\bar q_0}-\sqrt{\bar q_1})^2
-(n-1)p(\pi_0-\bar\pi_0)K\big\}.
\end{aligned}
\]
We can take the constant $c_{L,\varepsilon}>0$ large enough that
\begin{equation}\label{eq:joint_ub_pi}
\bar \pi_0 \leq \frac{\varepsilon}{4L(1+\varepsilon)}.
\end{equation}
Then we have 
\begin{equation}\label{eq:joint_H_lb}
\begin{aligned}
(\sqrt{\bar q_0} - \sqrt{\bar q_1})^2
&= \frac{L\bar\pi_0^2}{(\sqrt{1+(L-1)\bar \pi_0} + \sqrt{1-\bar \pi_0})^2}\\
& \geq \frac{L\bar \pi_0^2}{\big(2 + \tfrac{L-1}{2}\bar \pi_0\big)^2}\\
& \geq \frac{L}{4} \bar \pi_0^2 \Big(1-\frac{L-1}{2}\bar\pi_0\Big)
\end{aligned}
\end{equation}
where the first inequality comes from $\sqrt{1+x} \leq 1 + \tfrac{x}{2}$ and $\sqrt{1-\bar\pi_0}\leq 1$, and the second inequality comes from $(1+x)^{-2}\geq 1-2x$ for $x\geq 0$. Then combining \eqref{eq:joint_ub_pi} and \eqref{eq:joint_H_lb}, we obtain that
\[
(\sqrt{\bar q_0} - \sqrt{\bar q_1})^2
\geq \frac{L}{4}\bar \pi_0^2
\Big(1-\frac{\varepsilon}{4(1+\varepsilon)}\Big)
\]
which implies that 
\begin{equation}\label{eq:last_ineq}
(n-1)p(\sqrt{\bar q_0} - \sqrt{\bar q_1})^2
\geq (1+\varepsilon/2) \log n
\end{equation}
for sufficiently large $n$.

Since $\log(1+x) \leq x$, we have 
\[
\frac{1}{2}\log \frac{\bar q_0}{\bar q_1}
= \frac{1}{2}\log\Big(1 + \frac{L\bar \pi_0}{1-\bar \pi_0}\Big)
\leq \frac{L\bar\pi_0}{2(1-\bar\pi_0)}
\]
which implies that 
\begin{equation}\label{eq:third_ineq}
\frac{(1/2)\log(\bar q_0/\bar q_1)}
{(\sqrt{\bar q_0} - \sqrt{\bar q_1})^2}
\bar \pi_0
\leq
\frac{2}{(1-\bar\pi_0)
(1-\frac{L-1}{2}\bar\pi_0)}
<3.
\end{equation}
Moreover, from \eqref{eq:joint_ub_pi}, we have $t\leq1$, which implies that $K\geq t/2$.

Choosing a small constant $c_\varepsilon>0$ such that
$c_\varepsilon<1/2$ and
$(1-3c_\varepsilon)(1+\varepsilon/2)>1+\varepsilon/4$,
and taking $a = c_\varepsilon (n-1)p\pi_0$, we derive
\[
\begin{aligned}
&\mathbb{P}\big((u_k)_i \leq c_\varepsilon (n-1)p\pi_0\big)\\
&\leq \exp\Big\{(n-1)p\Big(
c_\varepsilon t\pi_0
-(\sqrt{\bar q_0}-\sqrt{\bar q_1})^2
-(\pi_0-\bar\pi_0)K
\Big)\Big\}\\
&=\exp\Big\{(n-1)p\Big(
c_\varepsilon t\bar\pi_0
-(\sqrt{\bar q_0}-\sqrt{\bar q_1})^2
+\underset{\geq 0}{\underbrace{(\pi_0-\bar\pi_0)}}\underset{\leq 0}{\underbrace{(c_\varepsilon t-K)}}
\Big)\Big\}\\
&\leq\exp\Big\{(n-1)p\Big(
c_\varepsilon t\bar\pi_0
-(\sqrt{\bar q_0}-\sqrt{\bar q_1})^2
\Big)\Big\}\\
&\leq \exp\Big\{-(1-3c_\varepsilon)(n-1)p
(\sqrt{\bar q_0}-\sqrt{\bar q_1})^2\Big\}\\
&\leq n^{-(1+\varepsilon/4)}.
\end{aligned}
\]
where the second last inequality comes from \eqref{eq:third_ineq} and the last inequality comes from \eqref{eq:last_ineq}.
Thus, by union bound, we obtain with probability at least $1-(L-1)n^{-\varepsilon/4}$,
\[
\alpha_s /n \geq c_\varepsilon\frac{n-1}{n}p\pi_0
\geq \frac{c_\varepsilon}{2}p\pi_0.
\]
After adjusting the constant $c_\varepsilon$, the proof is finished.
\end{proof}

Next, we show that $\tfrac{\beta_s}{np\pi_0}\rightarrow 0$ in probability, as $n\rightarrow\infty$.

\begin{lem}\label{lem:joint_beta}
If $\pi_0$ satisfies \eqref{eq:joint_pi}, then $\beta_s/n = o_\mathbb{P}(p\pi_0)$.
\end{lem}

\begin{proof}
It is true that 
\[
\mathbb{E}(u_k) = (n-1)p\pi_0 \vone_n
\]
and 
\[
\begin{aligned}
\var(X^{(k)}_{ij}) &= \mathbb{E}((X^{(k)}_{ij})^2) - (\mathbb E (X^{(k)}_{ij}))^2\\
&\leq  \mathbb{E}((X^{(k)}_{ij})^2)= p(q_0 + q_1) \leq 2p
\end{aligned}
\]
Then we have
\[
\begin{aligned}
\mathbb E \|\Pi u_k \|_2^2= \mathbb E \|\Pi(u_k -\mathbb E(u_k))\|_2^2&\leq \|\Pi\|^2 \mathbb E   \|u_k -\mathbb E(u_k)\|_2^2\\
&\leq  \mathbb E   \|u_k -\mathbb E(u_k)\|_2^2\\
& = \sum_{i} \var((u_k)_i)\\
& \leq 2n(n-1)p
\end{aligned}
\]
Thus, we obtain
\[
\frac{1}{n^2}\mathbb E [\beta_s^2] \leq \frac{2(L-1)p}{n}
\]
which implies that
\[
\frac{1}{n^2p^2\pi_0^2}\mathbb E [\beta_s^2] \leq \frac{2(L-1)}{np\pi_0^2} = O\big(\frac{1}{\log n}\big)
\]
Thus, we have $\beta_s/n = o_\mathbb{P}(p\pi_0)$. 
\end{proof}

Before discussing the order of $\rho_s$, we need the following lemma.

\begin{lem}
\label{lem:random_matrix_norm}
Let $E\in\CC^{n\times n}$ be a random Hermitian matrix
with $E_{ii}=0$, where the entries
$\{E_{ij}:i<j\}$ are independent and satisfy
\[
\mathbb E[E_{ij}]=0,\qquad
\mathbb E[|E_{ij}|^2]\leq p,\qquad
|E_{ij}|\leq 2.
\]
Then
\[
\mathbb E\|E\|
=O\left(\sqrt{np}+\sqrt{\log n}\right).
\]
\end{lem}
\begin{proof}
See \Cref{append:proof_of_lem_norm}.
\end{proof}

Last, we prove that $\tfrac{\rho_s}{np\pi_0}$ also goes to $0$ in probability as $n\rightarrow\infty$.

\begin{lem}\label{lem:joint_rho}
If $\pi_0$ satisfies \eqref{eq:joint_pi}, then $\rho_s/n = o_\mathbb{P}(p\pi_0)$.
\end{lem}

\begin{proof}
Denote 
\begin{equation*}
z^{\odot k} :=
\begin{pmatrix}
\omega^{-ks_1}\\
\vdots\\
\omega^{-ks_n}
\end{pmatrix},
\qquad
\Pi_k:=I-\frac1n z^{\odot k}(z^{\odot k})^H.
\end{equation*}
which is a special case of \eqref{eq:zkstar} when $G = \ZZ_L$.
Since we have
\[
\begin{aligned}
(C_k)_{ij} = \sum_{t=0}^{L-1} \omega^{-kt}(F_t)_{ij} =& -\mathds{1}\{\{i,j\}\in \Omega\}\sum_{t=0}^{L-1} \omega^{-kt}\mathds{1}\{Y_{ij}=t\}\\
=& -\mathds{1}\{\{i,j\}\in \Omega\} \omega^{-kY_{ij}}
\end{aligned}
\]
it follows that
\[
\mathbb E[(C_k)_{ij}] = -p \omega^{-k d_{ij}} \Big(q_0 + q_1\sum_{q\neq 0} \omega^{-kq}\Big) = -p\pi_0 z^{\odot k}_i\overline{z^{\odot k}_j}
\]
Note that $\diag(C_k) = \vzero$. Then it is true that 
\[
\mathbb E [C_k] = -p\pi_0 \Big(z^{\odot k} (z^{\odot k})^H - I_n\Big)
\]
Next we have $\Pi_k(\mathbb E [C_k] )\Pi_k = p\pi_0 \Pi_k \succeq 0$. Let $E_k:= C_k - \mathbb{E}(C_k)$. Since we have
\[
\begin{aligned}
\lambda_{\min}(\Pi_kC_k\Pi_k) &=  \lambda_{\min}(\Pi_k(E_k+\mathbb{E}(C_k))\Pi_k)\\
&\geq  \lambda_{\min}(\Pi_kE_k\Pi_k)\\
&\geq - \|\Pi_kE_k\Pi_k\|
\end{aligned}
\]
it follows that
\[
\rho_s \leq \max_{k\neq 0} \|\Pi_kE_k\Pi_k\| \leq  \max_{k\neq 0} \|E_k\|
\]
For any $i\neq j$, we have
\[
\begin{aligned}
\mathbb E[|(E_k)_{ij}|^2]
&=\mathbb E[|(C_k)_{ij}|^2]
-|\mathbb E[(C_k)_{ij}]|^2\\
&=p-p^2\pi_0^2\leq p,
\end{aligned}
\]
and
\[
|(E_k)_{ij}|
\leq |(C_k)_{ij}|+|\mathbb E[(C_k)_{ij}]|
\leq 1+p\pi_0\leq 2.
\]
From  \Cref{lem:random_matrix_norm}, it is follows that,
\[
\mathbb{E} \|E_k\| = O(\sqrt{np} + \sqrt{\log n})
\]
Thus, we obtain
\[
\frac{\rho_s}{np\pi_0} = O_{\mathbb P}\Big(\frac{1}{\sqrt{\log n}}\Big) = o_{\mathbb P}(1). \qedhere
\]
\end{proof}

Then we are ready to prove \Cref{thm:exact_joint}.
\begin{proof}[Proof of \Cref{thm:exact_joint}]
From \Cref{lem:joint_alpha}, \Cref{lem:joint_beta} and \Cref{lem:joint_rho}, we have with high probability,
\[
\alpha_s > c_\varepsilon n p \pi_0 > \rho_s + 2\sqrt{L-1}\beta_s
\]
that is, \eqref{eq:exact_condition} holds. Thus, SDP \eqref{eq:sdp} is tight with high probability.
\end{proof}

\subsection{Minimax Optimality} To conclude this section, we show that when $\pi_0$ falls below the threshold established
above, exact recovery of the true labels, up to a global cyclic
shift, is impossible in the minimax sense.

\begin{thm}[{\cite[Theorem~3.5]{chen2018projected}}]
Fix $L \in \ZZ_+$ and $\varepsilon \in (0,1)$, suppose that $p>c_1\tfrac{\log n}{n}$. If $\pi_0$ satisfies
\begin{equation}
\pi_0 < 2\sqrt{\frac{(1-\varepsilon) \log n}{Lpn}}
\end{equation}
then no method can recover the true solution up to a global shift.
\end{thm}

\section{Iterative minimal-rank rounding}
\label{sec:rounding}

In this section, we discuss how to round an optimal SDP solution
\(Z\in\cZ_G^n\setminus\cQ_G^n\). Suppose that the SDP returns a solution outside $\cQ_G^n$,
\(\Z^{(0)}=(Z_\pi^{(0)})_{\pi\in\widehat G\setminus\{\pi_0\}}\).
We define \(\|\Z\|_F^2:=\sum_{\pi\neq\pi_0}d_\pi\|Z_\pi\|_F^2\) and 
\[
\Z\bullet\Y:=\sum_{\pi\neq\pi_0}d_\pi\langle Z_\pi,Y_\pi\rangle.
\]
A natural rounding problem is to find the nearest point in $\cQ^n_G$,
\[
\min_{\Y \in\cQ_G^n}
\|\Z^{(0)}-\Y\|_F^2.
\]
Since every \(\Y\in\cQ_G^n\) satisfies
\(\|\Y\|_F^2=(|G|-1)n^2\) as shown below, this problem is equivalent to
\[
\max_{\Y\in\cQ_G^n}
\Z^{(0)}\bullet\Y
\]
Relaxing \(\cQ_G^n\) to the SDP feasible set motivates the iteration
\begin{equation}
\Z^{(r+1)}
\in
\argmax\{\Z^{(r)}\bullet\Z~\text{s.t.}~ \Z\in\cZ^n_G,\;Z_\pi\succeq0\ \forall\pi\neq\pi_0\}
\quad r=0,1,\dots.
\label{eq:rounding_iteration}
\end{equation}
Thus, each step can be interpreted as an SDP relaxation of the
nearest point projection problem onto $\cQ^n_G$. If an iterate returned by \eqref{eq:rounding_iteration} is of minimal rank, it is an exact nearest
minimal-rank point to the preceding iterate. 

A closely related fixed-point rounding method was proposed in ~\cite{felzenszwalb2022clustering,felzenszwalb2020iterated}. For the \texttt{Max-\(k\)-Cut} SDP
relaxation, the authors repeatedly solve SDPs obtained by relaxing a
closest-partition problem. Our procedure follows the same general
idea, but applies it to NUG-SDP feasible set, where the target discrete
solutions are characterized by rank-$d_\pi$ Fourier blocks.

\subsection{Fixed points and monotonicity.}

We now discuss the properties of the iterative rounding method. The following theorem shows that the Frobenius norm of the iterates is nondecreasing and bounded above, with the upper bound attained if and only if an iterate has minimal rank.

\begin{thm}
Let $\Z^{(0)}$ be an optimal solution of \eqref{eq:sdp} and let $\Z^{(r)}$ be the rounded solution obtained at the $r$-th iteration in \eqref{eq:rounding_iteration}. Then
\[
\|\Z^{(0)}\|_F^2\leq \|\Z^{(1)}\|_F^2 \leq \cdots\leq (|G|-1)n^2,
\]
where the upper bound is attained at the $r$-th iteration if and only if $Z_\pi^{(r)}$ has minimal rank, that is, $\rank(Z_\pi^{(r)})=d_\pi$ for every $\pi\in\widehat G\setminus\{\pi_0\}$.
\end{thm}

\begin{proof}

We first show that the Frobenius norm is nondecreasing along the iteration. Since \(\Z^{(r)}\) is feasible in the next subproblem, we have
\[
\Z^{(r)}
\bullet
\bigl(\Z^{(r+1)}-\Z^{(r)}\bigr)
\ge0,
\]
and therefore
\begin{equation*}
\|\Z^{(r+1)}\|_F^2 - \|\Z^{(r)}\|_F^2
\geq
\|\Z^{(r+1)}-\Z^{(r)}\|_F^2
\ge0.
\label{eq:rounding_monotonicity}
\end{equation*}

Then we show the upper bound. Given $\Z \in \cZ^n_G$ with $Z_\pi\succeq0$ for every $\pi\neq\pi_0$, for every feasible block, we can write
\[
Z_\pi=V_\pi V_\pi^H,
\qquad
V_\pi=
\begin{pmatrix}
V_{\pi,1}\\
\vdots\\
V_{\pi,n}
\end{pmatrix},
\]
where $V_{\pi,i} \in \CC^{d_\pi \times m}$ for some $m$. From the definition of $\cZ^n_G$, we have $(Z_\pi)_{ii} = V_{\pi,i} V_{\pi,i}^H = I_{d_\pi}$ and $\tr(Z_\pi) = nd_\pi$, then
\[
\begin{aligned}
\|Z_\pi\| = \|V_\pi V_\pi^H\| = \|V_\pi^H V_\pi\| &= \Big\|\sum_{j=1}^n V_{\pi,j}^H V_{\pi,j}\Big\|\\
&\leq \sum_{j=1}^n\| V_{\pi,j}\|^2 = n
\end{aligned}
\]
Let $\lambda_1\geq \lambda_2\dots,\geq\lambda_{nd_\pi} \geq 0$ denote the eigenvalues of $Z_\pi$. It follows that
\[
\begin{aligned}
\|Z_\pi\|_F^2 = \sum_j \lambda_j^2 \leq n\sum_j \lambda_j = n \tr(Z_\pi) = n^2 d_\pi
\end{aligned}
\]
where the inequality comes from $\lambda_j \leq \|Z_\pi\|\leq n$.
Thus, we derive
\[
\begin{aligned}
\|\Z\|_F^2 = \sum_{\pi \neq \pi_0} d_\pi \|Z_\pi\|_F^2 \leq  n^2\sum_{\pi \neq \pi_0} d_\pi^2 = n^2(|G|-1) 
\end{aligned}
\]
% \begin{equation*}\label{eq:rounding_inequality}
% Z_\pi\succeq0,
% \qquad
% (Z_\pi)_{ii}=I_{d_\pi},
% \qquad
% \tr(Z_\pi)=nd_\pi,
% \end{equation*}
% which implies that \(\|Z_\pi\|\le n\) and \(\|Z_\pi\|_F^2\le n^2d_\pi\), that is, $\|\Z\|_F^2\leq (|G|-1)n^2$.

Next, we show the upper bound is attained if and only if the solution has minimal rank.
If $\Z$ has minimal rank, then from \Cref{thm:rank-one}, we have $\Z \in \cQ^n_G$ which implies that $\|Z_\pi\|_F^2=n^2d_\pi$ and $\|\Z\|_F^2 = (|G|-1)n^2$.

Suppose for every $\pi\in\widehat G\setminus\{\pi_0\}$, we have $\|Z_\pi\|_F^2=n^2d_\pi$. Since $\lambda_i\leq n$, we have
\[
n\tr(Z_\pi) - \|Z_\pi\|_F^2 = n\sum_{i=1}^{nd_\pi}\lambda_i - \sum_{i=1}^{nd_\pi}\lambda_i^2 = \sum_{i=1}^{nd_\pi}\lambda_i(n-\lambda_i) = 0
\]
which implies that $\lambda_i \in \{0,n\}$. Since we have $\tr(Z_\pi)=nd_\pi$, there are exactly $d_\pi$ eigenvalues are nonzero, that is, $\Z$ has minimal rank.
\end{proof}

\begin{cor}
If $\Z^*\in\cZ^n_G$ with $\rank(Z^*_\pi)=d_\pi,\forall \pi \neq \pi_0$, then $\Z^*$ is a fixed point of \eqref{eq:rounding_iteration}, that is,
\[
\Z^*\in\argmax\{ \Z^{*}\bullet\Z.~\text{s.t.}~\Z\in\cZ^n_G,\;Z_\pi\succeq0\ \forall\pi\neq\pi_0\}
\]
\end{cor}

\begin{proof}
It is straightforward from the monotonicity of $\Phi$ and the fact that equality in the above bound holds if and only if $\Z$ has minimal rank, that is, $\rank(Z_\pi)=d_\pi,\forall \pi \neq \pi_0$.
\end{proof}

The above results partially explain the rank-promoting behavior of \eqref{eq:rounding_iteration}.
Every minimal-rank feasible point is an absorbing fixed point, although higher-rank fixed points may also exist. The following example gives a higher-rank fixed point.

\begin{ex}
When $n=3$ and $G=\ZZ_2$, the only nontrivial irreducible
representation has dimension $d_\pi=1$. Consider the following feasible solution $Z^* \in \cZ^3_{\ZZ_2}$ with $\rank(Z^*)=2>d_\pi=1$,
\[
Z^* =
\begin{pmatrix}
1 & -1 & 0\\
-1 & 1 & 0\\
0 & 0 & 1
\end{pmatrix}
\]
It is a fixed point of \eqref{eq:rounding_iteration}, since the SDP feasible set is $\{Z: \diag(Z) = \vone_3, Z\succeq 0\}$ and for any feasible $Z$, we have
\[
Z \bullet Z^* = 3 - 2 Z_{12} \leq 5 = \|Z^*\|_F^2
\]
\end{ex}

Hence, monotonicity alone does not guarantee convergence to minimal rank. Nevertheless, in our experiments the iteration reaches a minimal rank fixed point in most
instances, suggesting that many higher-rank fixed points have small or unstable attraction regions.

\section{Experiments}\label{sec:experiments}

In this section, we evaluate the proposed SDP formulation on three representative problems involving optimization over cyclic groups: (1) Max-2-Lin, (2) Joint Assignment, and (3) Multireference Alignment. These problems capture different aspects of the general framework studied in this paper, ranging from combinatorial constraint satisfaction to statistical signal recovery. For each problem, we also compare spectral rounding (SDP-spec.) and iterative rounding (SDP-iter.) applied to the SDP solutions. 
For spectral rounding, we extract the leading eigenvector of the first nontrivial Fourier block and round its entries to the nearest $L$-th roots of unity. 
The numerical results show that iterative rounding consistently outperforms spectral rounding across all the problems considered.

\subsection{Max-2-Lin}

We consider random Max-2-Lin instances over $\mathbb Z_L$ on two classes of
graphs: Erd\H{o}s--R\'enyi graphs and planar graphs. For the
Erd\H{o}s--R\'enyi instances, each pair of vertices is connected independently
with probability $p_{\Omega}$. For the planar instances, we instead generate
the constraint graph so that it is planar. Conditional on the graph, for every
edge $(i,j)$, the offset $c_{ij}$ is sampled independently and uniformly from
$\mathbb Z_L$. Thus, while the two graph classes differ in their underlying
topology, the Max-2-Lin equations are generated according to the same
distribution in both cases. All equations are assigned unit weight.

Throughout the experiments, we fix the number of variables to $n=25$, the
group size to $L=11$, and generate $N=50$ independent instances for each graph
class. For the Erd\H{o}s--R\'enyi experiments, we use a fixed observation
probability $p_{\Omega}\in\{0.2, 0.4\}$. For the planar experiments, the graph is generated
subject to the planarity constraint, and an independent uniformly distributed
offset $c_{ij}\in\mathbb Z_L$ is assigned to each edge.

We compare  NUG-SDP-based methods with the following benchmarks:

\begin{itemize}
    \item multistart local search;
    \item the algorithm of \cite{li_et_al:LIPIcs.ESA.2019.71}, denoted by LSZ;
    \item the algorithm of
    \cite{levialev_et_al:LIPIcs.APPROX-RANDOM.2017.18}, denoted by AL,
    which is applicable only to planar graphs.
\end{itemize}

\begin{figure}[h!]
\centering  
\subfigure[Square Grid]{
\includegraphics[width=0.45\linewidth]{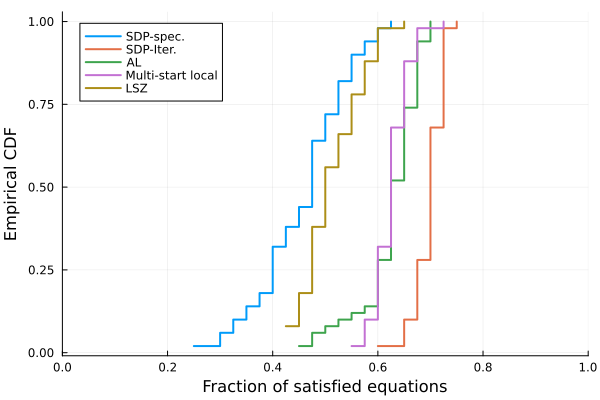}}
\subfigure[Triangulated Grid]{
\includegraphics[width=0.45\linewidth]{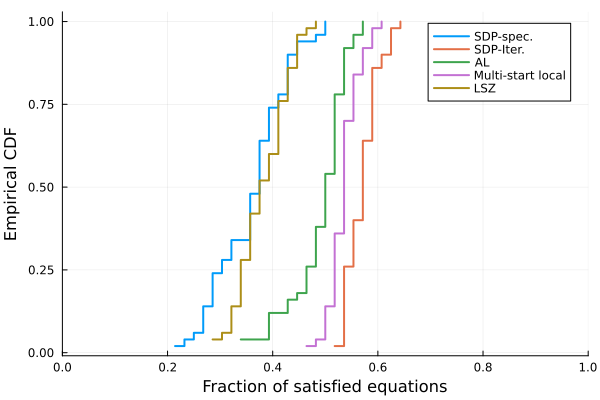}}

\caption{Max-2-Lin: empirical CDF of the fraction of satisfied equations on planar graphs.}
\label{fig:max-2-lin-planar}
\end{figure}

\begin{figure}[h!]
\centering  
\subfigure[Random graph: $p_\Omega=0.2$]{
\includegraphics[width=0.45\linewidth]{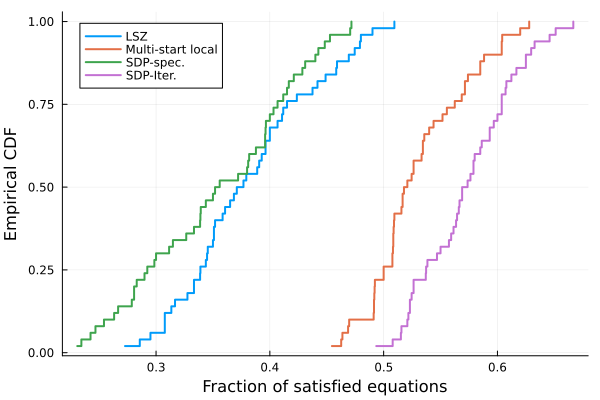}}
\subfigure[Random graph: $p_\Omega=0.4$]{
\includegraphics[width=0.45\linewidth]{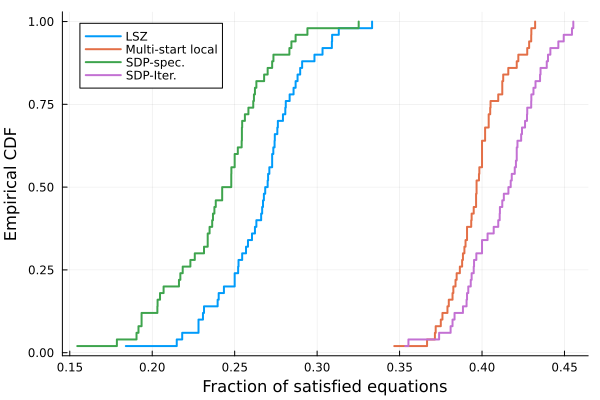}}

\caption{Max-2-Lin: empirical CDF of the fraction of satisfied equations on random graphs.}
\label{fig:max-2-lin-random}
\end{figure}

For each method, we evaluate the fraction of satisfied equations, which
coincides with the normalized objective value under unit weights. We summarize
the distribution of this quantity over the $N$ independent instances using
empirical cumulative distribution function (ECDF) plots. As shown in
Figures~\ref{fig:max-2-lin-planar} and \ref{fig:max-2-lin-random}, the SDP with iterative rounding achieves the best overall performance among the methods considered.

\subsection{Joint Assignment}

For the joint assignment problem, we first sample the ground-truth labels
$s_1^*,\ldots,s_n^*$ independently and uniformly from $\mathbb Z_L$.
Each pair of vertices is observed independently with probability
$p_{\Omega}$. For every observed pair $(i,j)$, the measurement equals the true
relative label
\[
s_i^* - s_j^* \bmod L
\]
with probability $\pi_0$; otherwise, it is sampled uniformly from
$\mathbb Z_L$. The pairwise costs are then defined by the negative
log-likelihood associated with this observation model. For each parameter
setting, we generate $N$ independent problem instances.

Throughout the experiments, we fix $n=30$, $L=15$, and $p_{\Omega}\in \{0.3,0.6\}$, and vary
$\pi_0$ to control the reliability of the observed pairwise measurements.
For each value of $\pi_0$, we generate $N=50$ independent instances and evaluate
all methods on the same collection of instances.

We compare NUG-SDP-based methods with the following benchmarks:
\begin{itemize}
    \item multistart local search;
    \item the projected power method (PPM) of \cite{chen2018projected};
    \item CEMP of \cite{lerman2022robust}.
\end{itemize}

\begin{figure}[h!]
	\centering  
        \subfigure[$p_\Omega = 0.3$]{
		\includegraphics[width=0.45\linewidth]{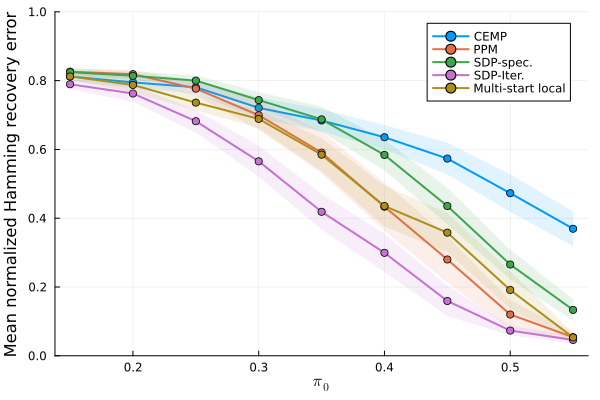}}
	\subfigure[$p_\Omega = 0.6$]{
		\includegraphics[width=0.45\linewidth]{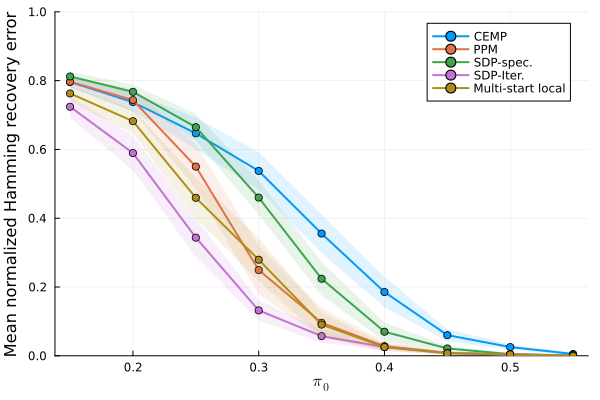}}
        \caption{Joint assignment: mean recovery error}
        \label{fig:JA_error}
\end{figure}

\begin{figure}[h!]
	\centering  
        \subfigure[$p_\Omega = 0.3, \pi_0 = 0.5$]{
		\includegraphics[width=0.45\linewidth]{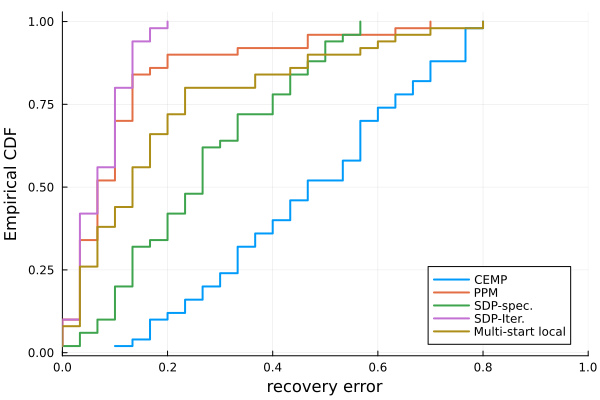}}
	\subfigure[$p_\Omega = 0.6, \pi_0 = 0.3$]{
		\includegraphics[width=0.45\linewidth]{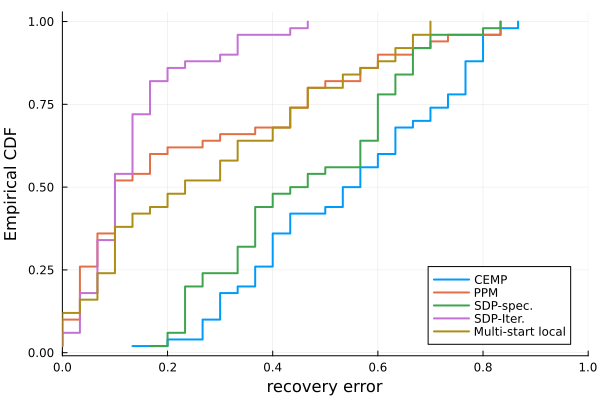}}
        \caption{Joint assignment: empirical CDF of recovery error.}
        \label{fig:JA_ecdf}
\end{figure}

For each method, we measure the recovery error by the normalized Hamming error up to a global shift, defined as
\[
\operatorname{Error}(\hat s,s^*)
=\min_{q\in\ZZ_L}\frac1n\sum_{i=1}^n
\mathds{1}\{\hat s_i+q\neq s_i^*\}.
\]
As shown in Figures~\ref{fig:JA_error} and \ref{fig:JA_ecdf}, the SDP with iterative rounding achieves
the lowest recovery error and exhibits the best overall performance among the methods considered.

\subsection{Multireference alignment}

\subsubsection{Synthetic Data}

We consider a random signal $x\in\RR^{m\times L}$, whose entries are sampled
independently from a standard Gaussian distribution. We generate observations
according to
\[
y_i = g_i x + \sigma \epsilon_i,
\qquad i\in[n],
\]
where $g_i$ acts by cyclically shifting the columns of $x$, and the entries of
$\epsilon_i$ are sampled independently from $\mathcal N(0,1)$. 
For each observation, we sample the group element $g_i$ independently and uniformly from $\ZZ_L$.
All $m$ channels of a given observation are shifted by the same group element $g_i$.  We define the signal-to-noise ratio (SNR) as
\(\mathrm{SNR}:=\|x\|_F^2/(mL\sigma^2)\).
For each parameter setting, we generate $N$ independent problem instances.

Throughout the experiments, we fix $m=5$ and $L=15$. We vary the noise level $\sigma$ or vary the sample size $n$ to control the difficulty of the recovery problem. For each value of $\sigma$ and $n$, we generate $N=50$ independent instances and evaluate all methods on the same collection of observations.

We compare NUG-SDP-based methods with the following benchmarks:
\begin{itemize}
    \item the expectation-maximization (EM) algorithm;
    \item the iterative method of \cite{gao2022iterative}, denoted by GZ.
    \item the synchronization-EM method \cite{janco2022accelerated}.
\end{itemize}

For each method, we report the signal reconstruction error after optimizing
over global cyclic shift. Specifically, if $\widehat x$
denotes the estimated signal, the recovery error is measured by
$\min_{q\in\mathbb Z_L}
\tfrac{\|g_q\widehat x-x\|_F}{\|x\|_F}$.

\begin{figure}[h!]
	\centering  
        \subfigure[SNR = $1/9$]{
		\includegraphics[width=0.45\linewidth]{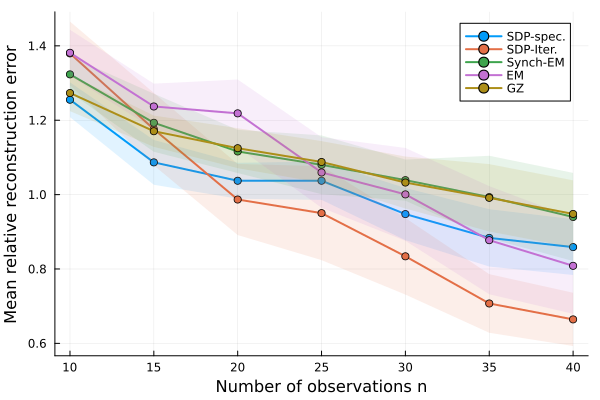}}
    	\subfigure[$n=30$]{
		\includegraphics[width=0.45\linewidth]{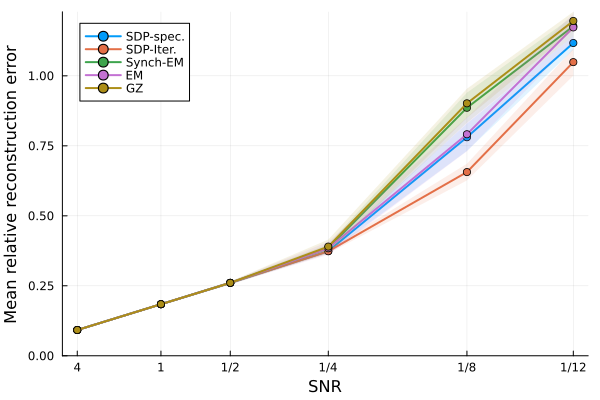}}
        
        \caption{MRA: mean recovery error.}
        \label{fig:mra_error}
\end{figure}

\begin{figure}[h!]
	\centering  
        \subfigure[$n=30$, SNR = $1/8$]{
		\includegraphics[width=0.45\linewidth]{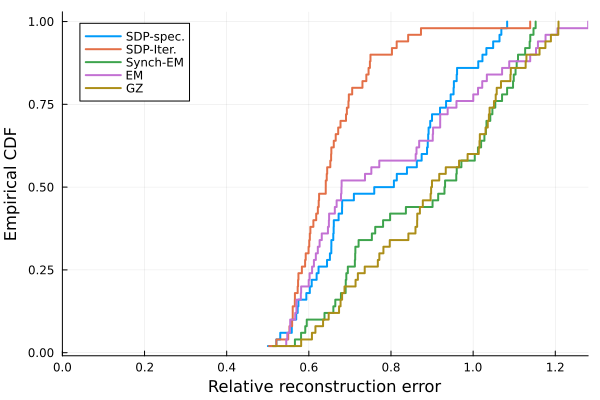}}
        \subfigure[$n=30$, SNR = $1/12$]{
		\includegraphics[width=0.45\linewidth]{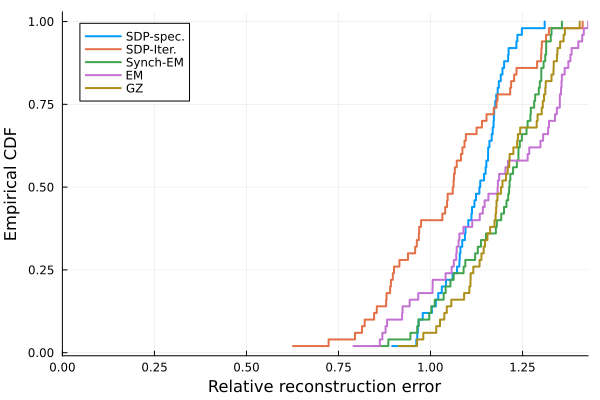}}
        
        \caption{MRA: empirical CDF of recovery error.}
        \label{fig:mra_ecdf}
\end{figure}

As shown in Figures~\ref{fig:mra_error} and \ref{fig:mra_ecdf}, the SDP with iterative rounding achieves the best
overall recovery performance among the methods considered.

\subsubsection{Image Alignment}

We next apply NUG-SDP to an image alignment problem. Let
$I_0\in\RR^{r\times r}$ denote a ground-truth image, and let $R_\theta$ denote the operator that rotates an image counterclockwise by an angle $\theta$. 
Throughout the image alignment experiments, we discretize the rotation group into $L$ possible rotations. For each observation, we sample $t$ uniformly from $\{0,\ldots,L-1\}$ and generate $ Y = R_{2\pi t/L} I_0 + \varepsilon$, where $\varepsilon\in\RR^{r\times r}$ has independent entries distributed as $\mathcal N(0,\sigma^2)$. Given $n$ independently generated noisy and rotated
images, our goal is to estimate the rotation associated with each observation, align the images accordingly, and recover the underlying image $I_0$.

For a recovered image $\widehat I$, we measure the reconstruction error up to a global rotation.  The noise level is characterized by the signal-to-noise ratio (SNR),
defined as $\mathrm{SNR}:=
\tfrac{\mathbb E(\mathrm{Signal}^2)}
     {\mathbb E(\mathrm{Noise}^2)}
=\tfrac{\|I_0\|_F^2}{r^2\sigma^2}$.
To formulate the image alignment problem within NUG-SDP framework over $\ZZ_L$, we
transform each image from Cartesian coordinates to polar coordinates. Under
this representation, each angular coordinate corresponds to a column of the
transformed image, while the entries within that column correspond to different
radial locations. Consequently, a planar rotation of the original image becomes
a cyclic shift along the angular coordinate. The image alignment problem can
therefore be treated as a multireference alignment problem over the cyclic
group $\mathbb Z_L$.

\begin{figure}[h!]
	\centering  
	% \subfigtopskip=1pt 
	% \subfigbottomskip=1pt
	% \subfigcapskip=-2pt 
        \subfigure[Image 1]{
		\includegraphics[width=0.3\linewidth]{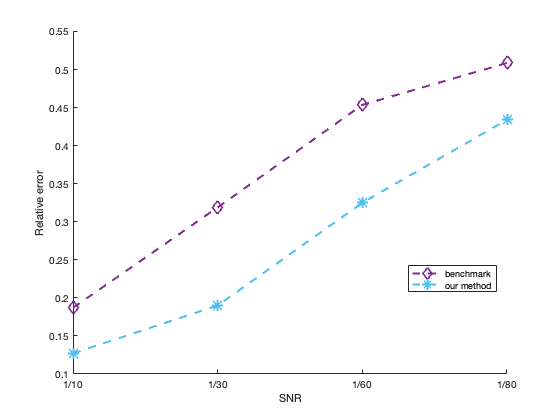}}
	\subfigure[Image 2]{
		\includegraphics[width=0.3\linewidth]{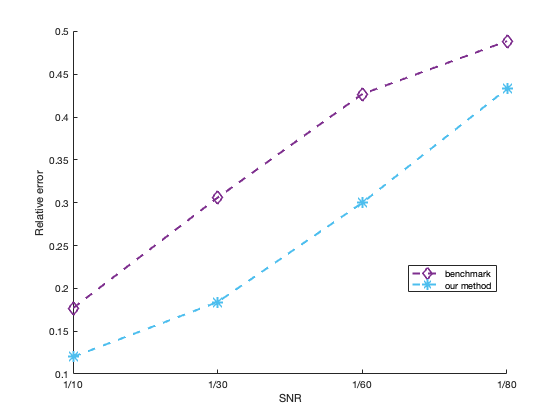}}
        \subfigure[Image 3]{
		\includegraphics[width=0.3\linewidth]{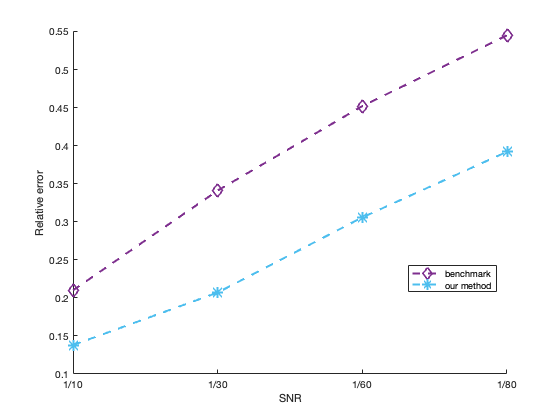}}
        \caption{Comparison of relative errors.}
        \label{fig:compare_relerr_img}
\end{figure}

As ground-truth images, we use random projections of the \emph{E.~coli} 70S
ribosome volume, each of size $129\times129$ ($r=129$). We consider SNR levels
ranging from $1/10$ to $1/80$, and fix $(n,L) = (40, 30)$. 
Our method consists of two steps. In Step 1, we implement NUG-SDP with iterative rounding to recover the rotation angles. In Step 2, we rotate the images back using the estimated angles, take their average, and then apply a mean filter for denoising. We compare our method with the bispectrum-based method \cite{ma2019heterogeneous}. The relative reconstruction errors are reported in Figure~\ref{fig:compare_relerr_img}.
Representative recovered images obtained by NUG-SDP with iterative rounding and the benchmark method are shown in Figures~\ref{fig:rec_2}, \ref{fig:rec_4} and
\ref{fig:rec_9}.

\captionof{figure}{Recovered image 1 at different SNR levels. For the high SNR case (first row), NUG-SDP recovers the rotation angles correctly for all samples. The ground-truth image is shown on the right.}

\begin{minipage}{0.55\linewidth}
    \centering
    \includegraphics[width=\linewidth]{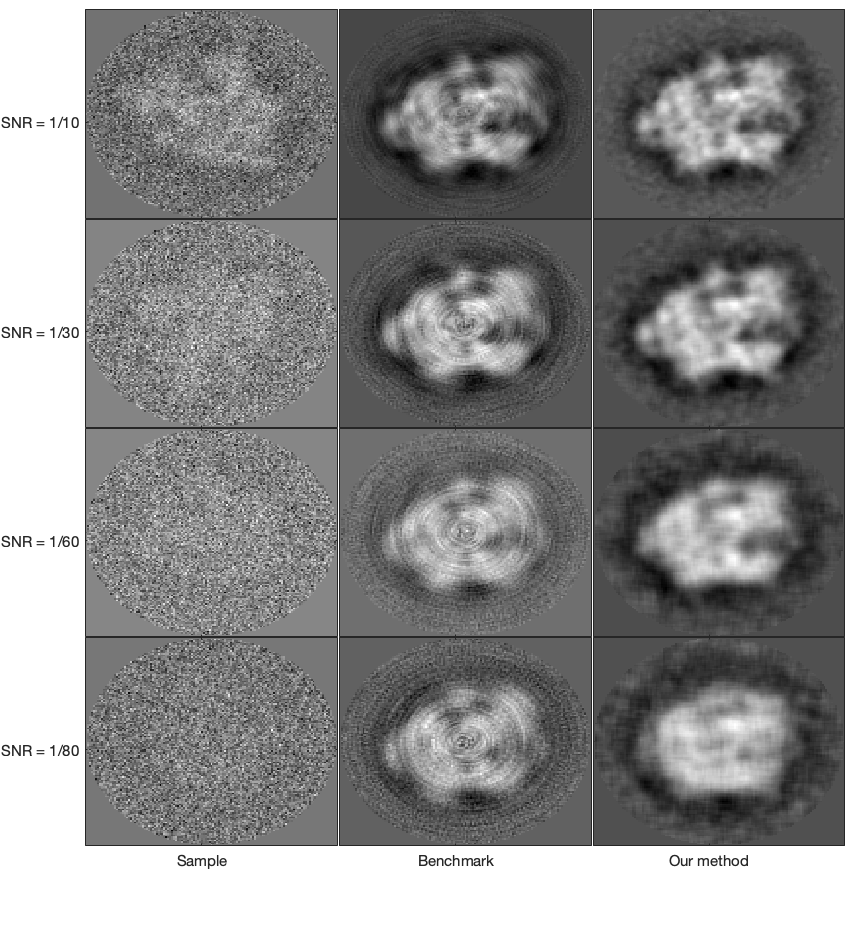}
    \label{fig:rec_2}
\end{minipage}
%\hspace{0.001\textwidth}
\begin{minipage}{0.3\linewidth}
    \centering
    \includegraphics[width=0.8\linewidth]{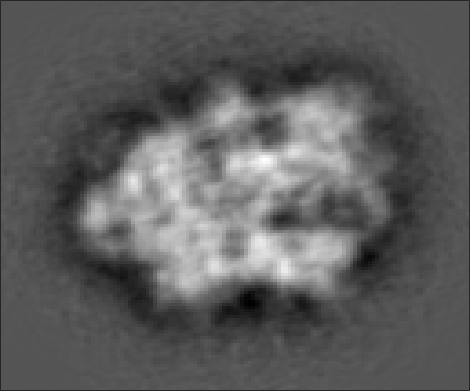}
    \label{fig:true_img_2}
\end{minipage}

\captionof{figure}{Recovered image 2 at different SNR levels. For the high SNR case (first row), NUG-SDP recovers the rotation angles correctly for all samples. The ground-truth image is shown on the right.}

\begin{minipage}{0.55\linewidth}
    \centering
    \includegraphics[width=\linewidth]{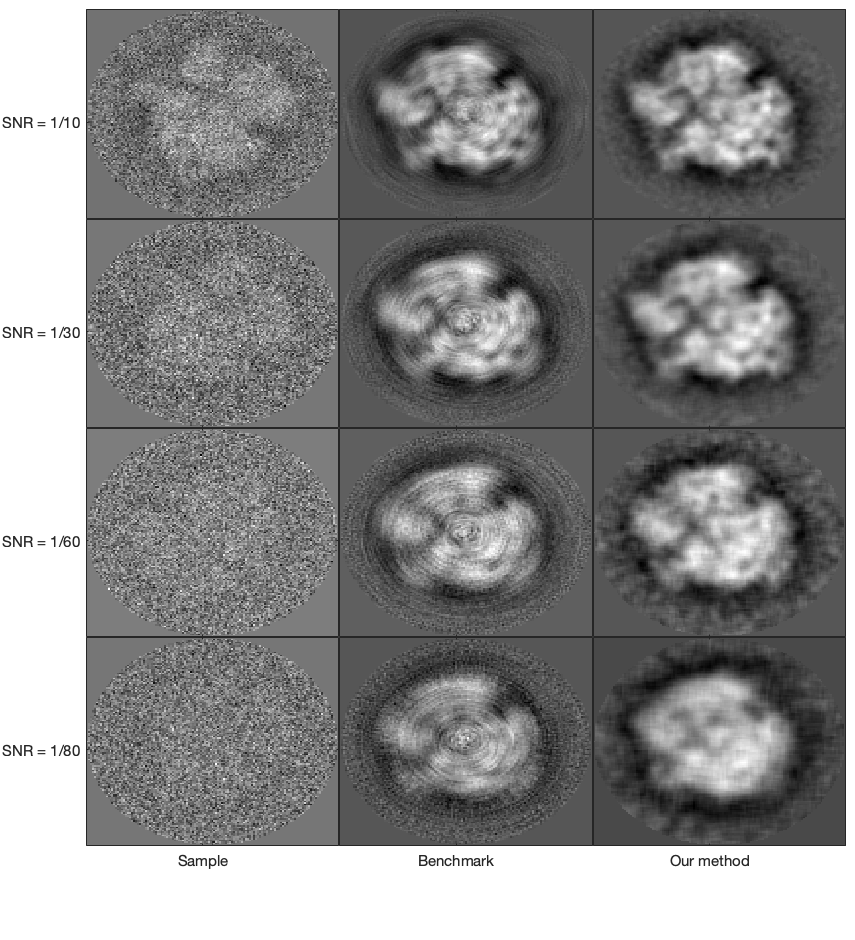}
    \label{fig:rec_4}
\end{minipage}
%\hspace{0.001\textwidth}
\begin{minipage}{0.3\linewidth}
    \centering
    \includegraphics[width=0.8\linewidth]{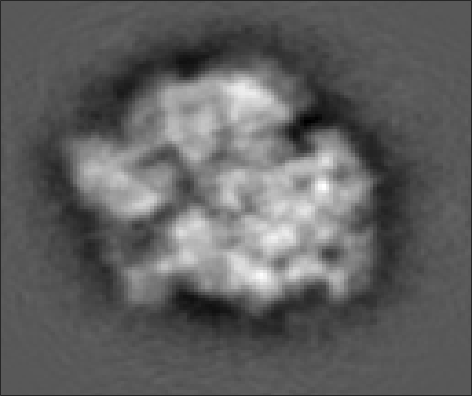}
    % \captionof{figure}{True image 2.}
    % right side
    \label{fig:true_img_4}
\end{minipage}

\captionof{figure}{Recovered image 3 at different SNR levels. For the high SNR case (first row), NUG-SDP recovers the rotation angles correctly for all samples. The ground-truth image is shown on the right.}

\begin{minipage}{0.55\linewidth}
    \centering
    \includegraphics[width=\linewidth]{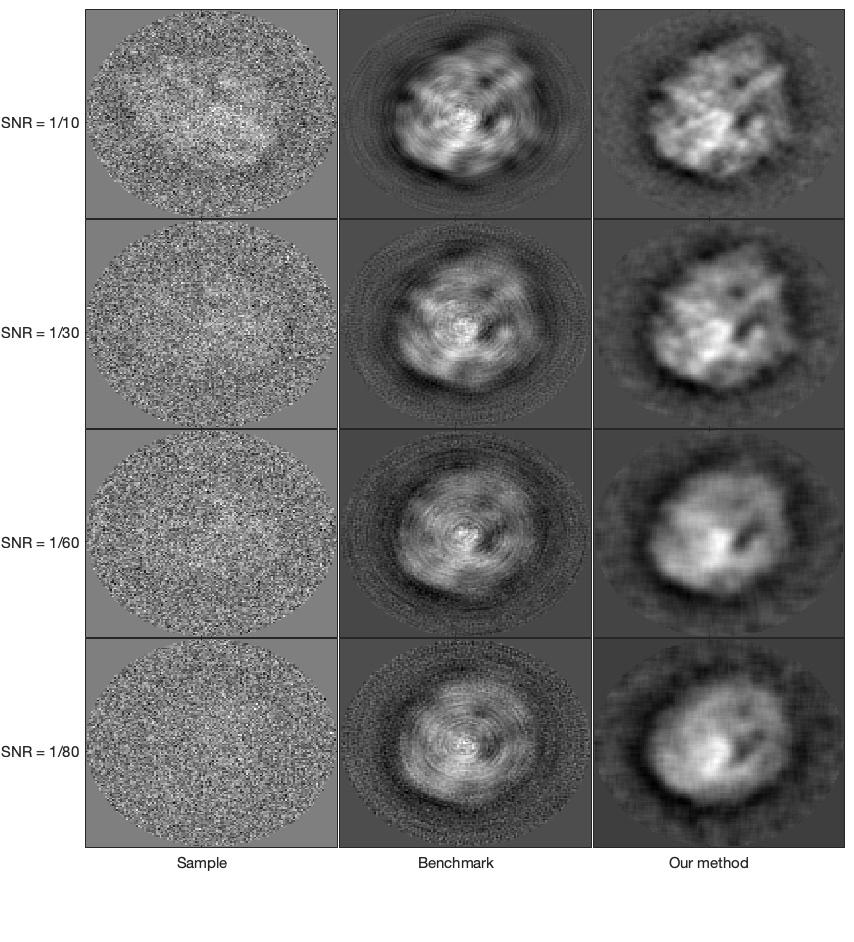}
    \label{fig:rec_9}
\end{minipage}
%\hspace{0.001\textwidth}
\begin{minipage}{0.3\linewidth}
    \centering
    \includegraphics[width=0.8\linewidth]{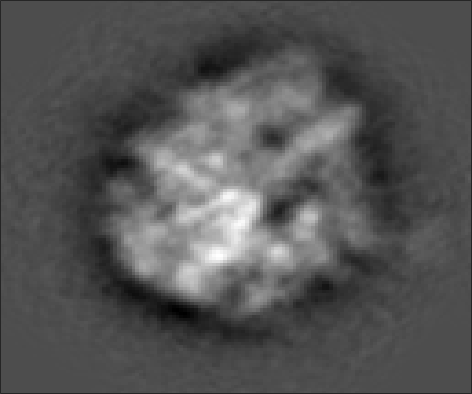}
    % \captionof{figure}{True image 3.}
    % right side
    \label{fig:true_img_9}
\end{minipage}

\section*{Acknowledgments}
The authors thank Xiaoming Huo for helpful discussions. The main mathematical analysis and initial results were developed independently by the authors. At an intermediate stage in the preparation of this manuscript, the authors used large language models to assist with language editing, exposition, checking selected mathematical derivations, and further developing and refining results based on the authors' original analysis. All AI-assisted content was carefully reviewed and verified by the authors. The authors assume full responsibility for all content.

% \section*{Declarations}

% \subsection*{Funding}
%  Dr.Cifuentes was partially supported by the Office of Naval Research, N00014-23-1-2631. Zhuorui Li was partially supported by ARC-ACO fellowship at Georgia Tech.

% \subsection*{Competing interests}
% The authors have no competing interests to declare that are relevant to the content of this article.

\bibliographystyle{abbrv}
\bibliography{ref.bib}

\begin{appendix}

\section{~}

\subsection{Background on Group Representations}
\label{append:group_representation}

We recall several standard results from finite-group representation theory used in our analysis. Let $G$ be a finite group with identity $e$, and let $\widehat G$ denote a complete set of inequivalent irreducible unitary representations of $G$. For each $\pi\in\widehat G$, let $d_\pi$ denote its dimension, and let $\pi_0$ denote the trivial representation.

\subsubsection{Plancherel's formula}\label{append:Plancherel}
For any function $f:G\to\CC$, its Fourier transform is defined by
\[
\hat f_\pi=\sum_{g\in G}f(g)\pi(g),
\qquad \pi\in\widehat G.
\]
Plancherel's formula over finite groups (see \cite[Corollary 4.3.5]{goodman2009symmetry}) states that for any $f,h:G\to\CC$,
\[
\sum_{g\in G}f(g)\overline{h(g)}
=\frac1{|G|}\sum_{\pi\in\widehat G}
d_\pi\operatorname{Tr}(\hat f_\pi\hat h_\pi^H).
\]
In particular, taking $h=f$, we obtain
\[
\sum_{g\in G}|f(g)|^2
=\frac1{|G|}\sum_{\pi\in\widehat G}
d_\pi\|\hat f_\pi\|_F^2.
\]
For the matrix-valued Fourier transform defined in this paper, applying the above identity entrywise gives
\[
\sum_{g\in G}\|M_g\|_F^2
=\frac1{|G|}\sum_{\pi\in\widehat G}
d_\pi\|\hat M_\pi\|_F^2,
\]
where $\hat M_\pi=\sum_{g\in G}M_g\otimes\pi(g)$.

Moreover, the dimensions of the irreducible representations satisfy
\[
\sum_{\pi\in\widehat G}d_\pi^2=|G|,
\]
see \cite[Corollary 4.4.5]{steinberg2012representation}.

\subsubsection{Schur's orthogonality relations}
\label{append:orthogonality}

By Schur's orthogonality relations
(see \cite[Theorem 4.2.8]{steinberg2012representation}),
for every nontrivial irreducible representation
$\pi\in\widehat G\setminus\{\pi_0\}$, we have
\[
\sum_{g\in G}\pi(g)=0_{d_\pi\times d_\pi}.
\]

\subsubsection{Irreducible representations of $G$ separate group elements}\label{append:irreducible}

The irreducible representations of $G$ separate group elements, i.e., for any $a,b\in G$,
\[
\pi(a)=\pi(b),\quad \forall\pi\in\widehat G
\quad\Longrightarrow\quad a=b.
\]
Indeed, let $\delta_a:G\to\CC$ denote the indicator function of $\{a\}$. Its Fourier transform satisfies $\hat\delta_{a,\pi}=\pi(a)$ for every $\pi\in\widehat G$. Thus, if $\pi(a)=\pi(b)$ for all $\pi\in\widehat G$, then $\delta_a$ and $\delta_b$ have the same Fourier transform. By the injectivity of the finite-group Fourier transform, we obtain $\delta_a=\delta_b$, and hence $a=b$ (see \cite[Proposition 5.5.5 and Remark 5.5.7]{steinberg2012representation}).

\subsection{Total number of PSD constraints} \label{append:less_PSD}

Suppose the potential functions have the undirected property, that is, \(f_{ij}(h) = f_{ji}(h^{-1})\), then \(C_\pi\) is Hermitian for any \(\pi\in\widehat G\). Denote by \(\bar\pi\) the complex-conjugate representation of \(\pi\). For each pair of inequivalent conjugate representations, we choose representatives such that \(\bar\pi(h)=\overline{\pi(h)}\).  
From the way \(\C\) is constructed in \eqref{eq:group_fourier}, we have
\[
C_\pi = \sum_{h\in G} F_h\otimes\pi(h)
= \sum_{h\in G} F_h\otimes\overline{\bar\pi(h)}
= \overline{C_{\bar\pi}}.
\]
which implies that \(C_\pi = C_{\bar\pi}^T\). For any \(\Z\in\cZ^n_G\), with \((\hat Z_h)_{h\in G}\in (\RR^{n\times n})^{|G|}\), we have
\[
Z_\pi = \sum_{h\in G} \hat Z_h\otimes\pi(h)
= \sum_{h\in G} \hat Z_h\otimes\overline{\bar\pi(h)}
= \overline{Z_{\bar\pi}} = Z_{\bar\pi}^T,
\]
where the last equality comes from \(Z_{\bar\pi}\in\mathbb{H}^{nd_\pi}\). Then it follows that
\[
d_\pi C_\pi \bullet Z_\pi
= d_\pi\tr(C_\pi Z_\pi)
= d_\pi\tr(Z_\pi^T C_\pi^T)
= d_{\bar\pi}C_{\bar\pi}\bullet Z_{\bar\pi}.
\]
Thus, in an implementation of \eqref{eq:sdp}, it only involves \((|\widehat G|-1+r)/2\) Hermitian PSD matrix variables, where $r$ is the number of nontrivial self-conjugate irreducible representations of $G$. 

If the potential functions do not satisfy the undirected property, we can obtain an equivalent formulation by symmetrizing each pair of potentials, thereby ensuring that the resulting functions satisfy this property.

\subsection{Proof of \Cref{thm:mra-converse}}\label{append_proof_of_mra_converse} 
We start with the following lemma.
\begin{lem}[{\cite[Theorem~6.1]{gao2022iterative}}]\label{lem:hypotheses}
Consider the following two hypotheses,
\[
H_0:y\sim \cN(x,\sigma^2 I),
\qquad
H_1:y\sim \cN(g_0^q x,\sigma^2 I),
\]
where $\|g_0^q x-x\|_F=\delta_{\min}$.
Under equal prior probabilities, the optimal testing
error is
\[
p_e=\mathbb P\left(
\cN(0,1)\geq\frac{\delta_{\min}}{2\sigma}
\right).
\]
\end{lem}

Next, we proceed to the proof. 
\begin{proof}[Proof of \Cref{thm:mra-converse}]
Following the proof of \cite[Theorem 6.1]{gao2022iterative},
we consider a restricted parameter space where $s_1^*=0$, and $s_i^*\in\{0,q\}$ for the remaining indices, where $q$ satisfies $\|g_0^q x-x\|_F=\delta_{\min}$ and $x$ is fixed. We say that $\hat s\equiv s^*$ up to a global shift if there exists $t\in\ZZ_L$ such that
\[
\hat s=s^*+t\vone_n\pmod L.
\]
Since $s_1^*=0$ in the restricted parameter space, for any estimator $\hat s$, let
\[
\tilde s=\hat s-\hat s_1\vone_n\pmod L.
\]
Then $\hat s\equiv s^*$ if and only if $\tilde s=s^*$. Thus, recovery up to a global shift is equivalent to exact recovery in this restricted parameter space.

For each $i\geq 2$, recovering $s_i^*$
is equivalent to testing the following two hypotheses,
\[
H_0:y_i\sim \cN(x,\sigma^2 I),
\qquad
H_1:y_i\sim \cN(g_0^q x,\sigma^2 I).
\]
Under equal prior probabilities, the optimal testing rule chooses the hypothesis whose mean is closer to $y_i$. From the assumption on $\sigma$, we have
\[
\frac{\delta_{\min}}{2\sigma}
\leq\frac{\sqrt{2\log n}}{1+\eta}.
\]
Thus, from the Gaussian tail lower bound and \Cref{lem:hypotheses}, we obtain the following lower bound for the optimal testing error $p_e$ for a single observation $y_i$,
\[
p_e=\mathbb P\left(
\cN(0,1)\geq\frac{\delta_{\min}}{2\sigma}
\right)\geq
\frac{c_\eta}{\sqrt{\log n}}
n^{-1/(1+\eta)^2}
\]
for some constant $c_\eta>0$ and sufficiently
large $n$.

Since the observations are independent, the optimal probability of recovering all these shifts is $(1-p_e)^{n-1}$.
Then we have
\[
\begin{aligned}
\inf_{\hat s}\sup_{s^*\in\ZZ_L^n}
\mathbb P(\hat s\not\equiv s^*)
&\geq 1-(1-p_e)^{n-1}\\
&\geq 1-\exp\big(-(n-1) p_e\big).
\end{aligned}
\]
Since $\eta>0$ is fixed, we have
\[
(n-1) p_e
\geq
\frac{c_\eta(n-1)}{\sqrt{\log n}}
n^{-1/(1+\eta)^2}
\longrightarrow\infty.
\]
Thus, we obtain
\[
\inf_{\hat s}\sup_{s^*\in\ZZ_L^n}
\mathbb P(\hat s\not\equiv s^*)
\longrightarrow1
\]
which implies that no method can recover the true solution.
\end{proof}

\subsection{Proof of \Cref{lem:random_matrix_norm}}\label{append:proof_of_lem_norm}

\begin{proof}
By \cite[Corollary 3.12 and Remark 3.13]
{bandeira2016sharp}, for a real symmetric random matrix
$X$ with independent centered entries above the diagonal,
bounded by a constant, we have
\[
\mathbb E\|X\|
=O\left(
\max_i\sqrt{\sum_j\mathbb E[X_{ij}^2]}
+\sqrt{\log n}
\right).
\]
To extend this result to Hermitian matrices, write
$E=A+\mathrm{i}(B-B^T)$, where
$A=\operatorname{Re}(E)$ and $B$ is the strictly
upper triangular part of $\operatorname{Im}(E)$.
The above bound applies to $A$ and to the following matrix
\[
\begin{pmatrix}
0 & B\\
B^T & 0
\end{pmatrix}.
\]
Since their entry variances are bounded by $p$,
we obtain
\[
\mathbb E\|A\|
=O\left(\sqrt{np}+\sqrt{\log n}\right),
\qquad
\mathbb E\|B\|
=O\left(\sqrt{np}+\sqrt{\log n}\right).
\]
The conclusion follows from
$\|E\|\leq\|A\|+2\|B\|$.
\end{proof}
\end{appendix}

\end{document}

%% file: notation.tex
\newcommand{\cD}{{\mathcal{D}}}

\newcommand{\cF}{{\mathcal{F}}}
\newcommand{\cG}{{\mathcal{G}}}

\newcommand{\cN}{{\mathcal{N}}}

\newcommand{\cQ}{{\mathcal{Q}}}

\newcommand{\cZ}{{\mathcal{Z}}}

\newcommand{\ZZ}{{\mathbb{Z}}}
\newcommand{\RR}{\mathbb{R}}

\newcommand{\vzero}{\mathbf{0}}
\newcommand{\vone}{{\mathbf{1}}}

\newcommand{\tr}{\text{Tr}}
\newcommand{\offdiag}{\text{offDiag}}

\newcommand{\vct}{\text{vec}}

\newcommand{\diag}{\text{diag}}

\newcommand{\argmax}{\text{argmax}}
\newcommand{\var}{\text{Var}}

\newcommand{\Z}{\boldsymbol{Z}}

\newcommand{\Y}{\boldsymbol{Y}}

\newcommand{\C}{\boldsymbol{C}}